\documentclass[a4paper,12pt]{amsart}
\usepackage{amssymb,amscd,amsmath,a4wide,graphicx,stmaryrd,fullpage,setspace,microtype,textcomp,tikz,accents,mathtools,abraces,tikz-cd}
\usepackage{hyperref}
\usepackage[shortlabels]{enumitem}
\definecolor{dark-red}{rgb}{0.5,0.15,0.15}
\hypersetup{
	colorlinks   = true, 
	urlcolor     = dark-red, 
	linkcolor    = black, 
	citecolor   = dark-red 
}
\usepackage[T1]{fontenc}
\usepackage{lmodern}
\usepackage[numbers,sort&compress]{natbib}

\title{Natural homotopy of multipointed d-spaces}

\author[P. Gaucher]{Philippe Gaucher}

\address{Universit\'e Paris Cit\'e, CNRS, IRIF, F-75013, Paris, France}

\urladdr{https://www.irif.fr/{\~{}}gaucher} 

\makeatletter
\@namedef{subjclassname@2020}{%
	\textup{2020} Mathematics Subject Classification}
\makeatother


\swapnumbers

\newcommand{\C}{\mathcal{C}}
\newcommand{\D}{\mathcal{D}}
\newcommand{\W}{\mathcal{W}}
\newcommand{\F}{\mathcal{F}}
\newcommand{\K}{\mathcal{K}}

\newcommand{\p}{\times}
\renewcommand{\vec}{\overrightarrow}
\renewcommand{\P}{\mathbb{P}}

\DeclareMathOperator{\diag}{Diag}
\DeclareMathOperator{\h}{\underline{\mathbf{h}}}

\newcommand{\cont}{{\vec{\mathrm{Sp}}}}
\newcommand{\opcont}{{\overleftarrow{\mathrm{Sp}}}}
\newcommand{\discont}{{\vec{\Omega}}}
\newcommand{\NT}{\overrightarrow{\mathrm{NT}}}
\newcommand{\dtop}{{\brm{Flow}}}
\newcommand{\dcat}{{\mathrm{cat}}}
\newcommand{\moore}{{\mathbb{M}}}
\newcommand{\lmoore}{\mathbb{M}_!}
\newcommand{\glob}{{\rm{Glob}}}

\newtheorem*{thmN}{Theorem}

\newtheorem{thm}{Theorem}[section]
\newtheorem{prop}[thm]{Proposition}
\newtheorem{lem}[thm]{Lemma}
\newtheorem{cor}[thm]{Corollary}
\newtheorem{defnot}[thm]{Definition and notation}
\newcommand{\bdn}{\begin{defnot}}
	\newcommand{\edn}{\end{defnot}}
\newcommand{\bp}{\begin{prop}}
	\newcommand{\ep}{\end{prop}}
\newcommand{\bth}{\begin{thm}}
	\renewcommand{\eth}{\end{thm}}
\newcommand{\bpf}{\begin{proof}}
	\newcommand{\epf}{\end{proof}}
\newcommand{\bc}{\begin{cor}}
	\newcommand{\ec}{\end{cor}}

\theoremstyle{definition}
\newtheorem{defn}[thm]{Definition}

\newcommand{\bd}{\begin{defn}}
	\newcommand{\ed}{\end{defn}}
\newtheorem{nota}[thm]{Notation}

\renewcommand{\top}{{\mathbf{Top}}}
\newcommand{\ho}{{\mathbf{Ho}}}
\newcommand{\iso}{\cong}
\newcommand{\vI}{\vec{I}}
\renewcommand{\leq}{\leqslant}
\renewcommand{\geq}{\geqslant}
\newcommand{\tr}[1]{{\langle{#1}\rangle}}
\DeclareMathOperator{\dt}{dt}
\newcommand{\ct}[1]{{\langle{#1}\rangle}}

\newcommand{\brm}[1]{\rm{\mathbf{#1}}}
\newcommand{\ttop}{{\brm{TOP}}}
\newcommand{\globM}{{\rm{Glob}}^{top}}

\DeclareMathOperator{\id}{Id}
\DeclareMathOperator{\Obj}{Obj}

\newcommand{\liminj}{\varinjlim}
\newcommand{\limproj}{\varprojlim}
\newcommand{\cat}{{\mathbf{Cat}}}
\newcommand{\rest}{\!\upharpoonright\!}
\DeclareMathOperator{\carrier}{Carrier}
\newcommand{\ptop}[1]{{\brm{{#1}dTop}}}

\newcommand{\cocartesian}{\arrow[lu, phantom, "\ulcorner"{font=\Large}, pos=0]}

\makeatletter
\def\varholim@#1#2{%
	\vtop{\m@th\ialign{##\cr
			\hfil$#1\operator@font holim$\hfil\cr
			\noalign{\nointerlineskip\kern1.5\ex@}#2\cr
			\noalign{\nointerlineskip\kern-\ex@}\cr}}%
}
\def\holimproj{%
	\mathop{\mathpalette\varholim@{\leftarrowfill@\textstyle}}\nmlimits@
}
\def\holiminj{%
	\mathop{\mathpalette\varholim@{\rightarrowfill@\textstyle}}\nmlimits@
}
\makeatother

\makeatletter
\def\@textbottom{\vskip \z@ \@plus 1pt}
\let\@texttop\relax
\makeatother

\makeatletter
\newcommand*{\@opargbegintheorem}[3]{\trivlist
	\item[\hskip \labelsep{\bfseries #1\ #2}] \textbf{(#3)}\ \itshape}
\makeatother

\newcommand{\ot}{\otimes}
\newcommand{\ttt}{two-out-of-three property}
\DeclareMathOperator{\natgl}{\underline{nat}^{gl}}

\allowdisplaybreaks

\begin{document}

\begin{abstract} 
	We identify Grandis' directed spaces as a full reflective subcategory of the category of multipointed $d$-spaces. When the multipointed $d$-space realizes a precubical set, its reflection coincides with the standard realization of the precubical set as a directed space.  The reflection enables us to extend the construction of the natural system of topological spaces in Baues-Wirsching's sense from directed spaces to multipointed $d$-spaces. In the case of a cellular multipointed $d$-space, there is a discrete version of this natural system which is proved to be bisimilar up to homotopy. We also prove that these constructions are invariant up to homotopy under globular subdivision. These results are the globular analogue of Dubut's results. Finally, we point the apparent incompatibility between the notion of bisimilar natural systems and the q-model structure of multipointed $d$-spaces and we give some suggestions for future works.
\end{abstract}

\maketitle
\tableofcontents
\hypersetup{linkcolor = dark-red}

\section{Introduction}

\subsection*{Presentation}

Directed spaces are topological spaces equipped with a set of continuous paths called directed paths closed under non-decreasing reparametrization and under composition and containing the constant paths \cite[Definition~1.1]{mg}. They are one of the geometric models of concurrency studied in Directed Algebraic Topology (DAT) \cite{DAT_book}. In his doctoral dissertation \cite{dubut_PhD}, Jérémy Dubut studies some natural systems in Baues-Wirsching's sense associated with a directed space. The purpose of this approach is to obtain global invariants of directed spaces, meaning algebraic objects capable of reducing the size of the state space without destroying the causal structure. Preserving this causal structure is the central problem of DAT, due to the non-conventional behavior of the directed segment. Indeed, the preservation of the causal structure requires extreme caution before contracting a segment in the direction of time. One of Dubut's results is the proof that some natural system associated with the realization as a directed space of a cubical complex is bisimilar to a discrete version of this construction \cite[Section~8.3.2 and 8.3.3]{dubut_PhD} (see also \cite[Theorem~1]{NaturalHomology}). An immediate consequence is the invariance under cubical subdivision of these constructions, up to bisimilarity. This result proves the relevance of the notion of bisimilarity of natural systems for DAT. Another interesting feature of this approach is the homological properties of the category of all small diagrams of e.g. abelian groups, which seems to pave the way towards calculations using generalizations of homological algebra.

The aim of this paper is to present the globular analogue of Dubut's results. The interest is twofold. It is a new step towards a unified formalism of the globular and cubical approaches of DAT. It also opens the way to computational techniques by rewriting a continuous directed path as a finite sequence of discrete globular cells. The starting point is the following observation.

\begin{thmN} (Theorem~\ref{thm:reflection})
	The directed spaces in the sense of \cite[Definition~1.1]{mg} are a full reflective subcategory of the category of multipointed $d$-spaces in the sense of \cite[Definition~6]{Moore3}.
\end{thmN}

The reflection $\cont$, constructed in Proposition~\ref{prop:cont}, takes a multipointed $d$-space to a directed space such that the directed paths are all finite compositions of pieces of execution paths. The image by the functor $\cont$ of the realization of a precubical set as a multipointed $d$-space is the standard realization of a precubical set as a directed space (see the end of Section~\ref{sec:from-multi-to-continuous}). This enables us to define a natural system in Baues-Wirsching's sense of topological spaces $\NT(X)$ associated with a multipointed $d$-space $X$ which extends the construction on directed spaces presented in \cite[Section~6.4]{dubut_PhD}. In the cellular case, we find a discrete version $\NT_d(X)$ of this natural system by working with the centers of the globular cells, exactly like in the cubical setting in which only the centers of the cubical cells are considered. We prove that $\NT_d(X)$ is bisimilar up to homotopy to $\NT(X)$. Unlike what happens in the cubical setting, the globular version of these results does not require to restrict to objects equipped with a global order. Finally, we also prove that these constructions are invariant up to homotopy under the T-homotopy equivalences as introduced in \cite[Definition~4.10]{diCW} and which are preferably called \textit{globular subdivisions} in this paper (see Figure~\ref{fig:globular-sdb}). These results are summarized in the following statements:

\begin{thmN} (Theorem~\ref{thm:open-up-to-homotopy})
	Let $X$ be a cellular multipointed $d$-space. There exists a map of natural systems of topological spaces from $\NT(X)$ to its discrete version $\NT_d(X)$ which is open up to homotopy. Consequently, these two natural systems are bisimilar up to homotopy.
\end{thmN}

\begin{thmN} (Theorem~\ref{thm:final})
	A globular subdivision is a map $f:X\to Y$ between cellular multipointed $d$-spaces inducing a homeomorphism between the underlying spaces (see Figure~\ref{fig:globular-sdb}). In this situation, we prove that the natural systems of topological spaces $\NT(X)$, $\NT(Y)$, $\NT_d(X)$ and $\NT_d(Y)$ are bisimilar up to homotopy.
\end{thmN}

The very last result of this paper is a negative result. Proposition~\ref{prop:counterexample} establishes that there exist two cellular multipointed $d$-spaces which are weakly equivalent in the q-model structure of multipointed $d$-spaces and such that the associated natural systems are not bisimilar up to homotopy. The counterexample does not use far-fetched multipointed $d$-spaces. They are both obtained by starting from a very simple finite loop-free precubical set of dimension $1$. This negative result means that the two approaches of DAT which are on one hand the model structures on multipointed $d$-spaces of \cite{QHMmodel} and on the other hand the notion of bisimilar natural systems seem to be incompatible. The paper is concluded by some explanations about this apparent incompatibility.

\subsection*{Outline of the paper}

Section~\ref{sec:Moore-regular} expounds some very basic notions and results about Moore paths and the regular ones. Some facts are generalizations or adaptations of statements coming from \cite{Moore2,Moore3}. It is important to notice that, in this paper, we do work with constant paths as well. 

Section~\ref{sec:from-multi-to-continuous} constructs the functor $\cont$ from the category of multipointed $d$-spaces to the category of directed spaces. In this paper, a directed space means a $d$-space in Grandis' sense \cite[Definition~1.1]{mg}. We prefer this terminology, also sometimes used in the literature, to avoid any confusion with multipointed $d$-spaces. In fact, we prove that the category of directed spaces is a full reflective subcategory of the category of multipointed $d$-spaces in Theorem~\ref{thm:reflection}, the functor $\cont$ being the reflection.

Section~\ref{sec:directed-paths-cellular-case} explores the basic properties of the directed paths of a cellular multipointed $d$-space. The important notion of discrete trace is introduced and it is proved that every directed path of a cellular multipointed $d$-space has a unique discrete trace. 

Section~\ref{sec:def-bisim-up-to-homotopy)} recalls some facts about the bisimilarity of diagrams and introduces the notions of bisimilarity up to homotopy and of open map up to homotopy. 

Section~\ref{sec:construction-NT} presents the construction of a functor from the category of multipointed $d$-spaces to the category of all small diagrams of topological spaces which takes a multipointed $d$-space $X$ to its natural system $\NT(X)$. 

Section~\ref{sec:path2trace} recalls an important fact about the passage from directed paths to traces which will be used in this paper. Proposition~\ref{prop:colim} is certainly not the most general valid statement. It is adapted to the use made in this paper. 

Section~\ref{sec:construction-NTd} introduces the discrete version $\NT_d(X)$ of $\NT(X)$ for a cellular multipointed $d$-space $X$ and establishes Theorem~\ref{thm:open-up-to-homotopy}. 

Section~\ref{sec:globular-sbd} is devoted to the globular subdivisions and to the proof of Theorem~\ref{thm:final}. 

Section~\ref{sec:counter-example} describes the incompatibility between the q-model structure of multipointed $d$-spaces and the notion of bisimilar natural systems. It is concluded by a discussion by studying the cases of two directed spaces not containing globes.

Section~\ref{sec:future} speculates about this incompatibility and gives some hints for future works.

Appendix~\ref{sec:opcont} is devoted to the proof that the category of directed spaces is also a full coreflective subcategory of the category of multipointed $d$-spaces in Theorem~\ref{thm:coreflection}. This fact is mentioned for the sake of completeness. It is not used in the core of the paper.

\subsection*{Acknowledgments}

I thank Jérémy Dubut for answering my questions about his doctoral dissertation. I thank the anonymous referee for the careful reading of the paper.

\section{Moore paths and regular paths}
\label{sec:Moore-regular}

We refer to \cite{TheBook} for locally presentable categories. The category $\top$ denotes either the category of \textit{$\Delta$-generated spaces} or of \textit{$\Delta$-Hausdorff $\Delta$-generated spaces} (cf. \cite[Section~2 and Appendix~B]{leftproperflow}). It is Cartesian closed by a result due to Dugger and Vogt recalled in \cite[Proposition~2.5]{mdtop} and locally presentable by \cite[Corollary~3.7]{FR}. The internal hom is denoted by $\ttop(-,-)$. The right adjoint of the inclusion from $\Delta$-generated spaces to general topological spaces is called the $\Delta$-kelleyfication. 

The notations $\ell,\ell',\ell_i,L,\dots$ mean a nonnegative real number. $[\ell,\ell']$ denotes a segment: unless specified, it is always understood that $\ell\leq\ell'$. Let $\mu_{\ell}:[0,\ell]\to [0,1]$ be the homeomorphism defined by $\mu_\ell(t) = t/\ell$ for $\ell>0$. 

Let $\mathcal{M}(\ell,\ell')$ be the set of non-decreasing surjective maps from $[0,\ell]$ to $[0,\ell']$ equipped with the $\Delta$-kelleyfication of the relative topology induced by the set inclusion $\mathcal{M}(\ell_1,\ell_2) \subset \ttop([0,\ell_1],[0,\ell_2])$. Note that all non-decreasing surjective maps between segments are continuous. By an argument similar to the one of \cite[Proposition~2.5]{Moore2}, we see that the topology of $\mathcal{M}(\ell,\ell')$ coincides with the compact-open topology and with the pointwise convergence topology.

\begin{nota} \label{nota:nondecreasing-surjective-ot}
	Let $\phi_i\in \mathcal{M}(\ell_i,\ell'_i)$ for $n\geq 1$ and $1\leq i \leq n$. Then the map
\[
\phi_1 \ot \dots \ot \phi_n : \big[0,\sum_i \ell_i\big] \longrightarrow \big[0,\sum_i \ell'_i\big]
\]
denotes the non-decreasing surjective map defined by 
\[
(\phi_1 \ot \dots \ot \phi_n)(t) = 
\begin{cases}
		\phi_1(t) & \hbox{if } 0\leq t\leq \ell_1\\
		\phi_2(t-\ell_1)+\ell'_1 & \hbox{if } \ell_1\leq t\leq \ell_1+\ell_2\\
		\dots \\
		\phi_i(t-\sum_{j<i}\ell_j) + \sum_{j<i}\ell'_j& \hbox{if } \sum_{j<i}\ell_j\leq t \leq \sum_{j\leq i}\ell_j\\
		\dots\\
		\phi_n(t-\sum_{j<n}\ell_j) + \sum_{j<n}\ell'_j & \hbox{if } \sum_{j<n}\ell_j\leq t \leq \sum_{j\leq n}\ell_j.
\end{cases}
\] 
\end{nota}

\bd \label{def_Moore_path}
Let $U$ be a topological space. A \textit{(Moore) path} of $U$ consists of a continuous map $\gamma:[0,\ell]\to U$. The real number $\ell$ is called the \textit{length} of $\gamma$. Remember that in this paper, and unlike some previous papers like \cite{Moore2,Moore3,RegularMoore}, we also consider Moore paths of length $0$.
\ed

\bd Let $\gamma_1:[0,\ell_1]\to U$ and $\gamma_2:[0,\ell_2]\to U$ be two Moore paths of a topological space $U$ such that $\gamma_1(\ell_1)=\gamma_2(0)$. The \textit{Moore composition} $\gamma_1*\gamma_2:[0,\ell_1+\ell_2]\to U$ is the Moore path defined by 
\[
(\gamma_1*\gamma_2)(t)=
\begin{cases}
	\gamma_1(t) & \hbox{ for } t\in [0,\ell_1]\\
	\gamma_2(t-\ell_1) &\hbox{ for }t\in [\ell_1,\ell_1+\ell_2].
\end{cases}
\]
The Moore composition of Moore paths is strictly associative.
\ed

\bd
Let $\gamma_1$ and $\gamma_2$ be two continuous maps from $[0,1]$ to some topological space such that $\gamma_1(1)=\gamma_2(0)$. The composite defined by 
\[
(\gamma_1 *_N \gamma_2)(t) = \begin{cases}
	\gamma_1(2t)& \hbox{ if }0\leq t\leq \frac{1}{2},\\
	\gamma_2(2t-1)& \hbox{ if }\frac{1}{2}\leq t\leq 1
\end{cases}
\]
is called the \textit{normalized composition}. One has \[\gamma_1 *_N \gamma_2 = (\gamma_1\mu_{1/2}) * (\gamma_2\mu_{1/2}).\] 
\ed

\bp \label{lem-1} 
Let $\ell_1,\ell'_1,\dots,\ell_n,\ell'_n\geq 0$. Let $U$ be a topological space. Let $\gamma_i:[0,\ell'_i]\to U$ be $n$ continuous maps with $1\leq i \leq n$ and $n\geq 1$. Let $\phi_i:[0,\ell_i]\to [0,\ell'_i]$ be a non-decreasing surjective map for $1\leq i \leq n$. Then we have 
\[
\left(\gamma_1*\dots *\gamma_n\right)(\phi_1\ot\dots\ot\phi_n) = (\gamma_1\phi_1)*\dots *(\gamma_n\phi_n).
\]
\ep

\bpf 
For $\sum_{j<i}\ell_j\leq t \leq \sum_{j\leq i}\ell_j$, we have 
\[\begin{aligned}
	\big(\gamma_1*\dots *\gamma_n\big)(\phi_1 \ot \dots \phi_n)(t) &= \gamma_i\big((\phi_1 \ot \dots \ot \phi_n)(t)- \sum_{j<i}\ell'_j\big)\\
	&= \gamma_i\big(\big(\phi_i\bigl(t-\sum_{j<i}\ell_j\big) + \sum_{j<i}\ell'_j\big)- \sum_{j<i}\ell'_j\big)\\
	&= \gamma_i\big(\phi_i\bigl(t-\sum_{j<i}\ell_j\big)\big)\\
	& = \big((\gamma_1\phi_1)*\dots *(\gamma_n\phi_n)\big)(t),
\end{aligned}\]
the first and fourth equalities by definition of the Moore composition, the second equality by definition of $\phi_1 \ot \dots \phi_n$, and the third equality by algebraic simplification. 
\epf

We need in Theorem~\ref{continuous} a generalization of Proposition~\ref{lem-1}. We must at first extend the notion of tensor product of non-decreasing surjective maps between segments to \textit{composable} non-decreasing continuous maps between segments as follows (we treat only the case of two continuous maps, which is the only one used in this paper).

\begin{nota} \label{nota:nondecreasing-continuous-ot}
	Let $\ell_1,\ell'_1,\ell_2,\ell'_2\geq 0$. Consider two non-decreasing continuous maps $\phi_i:[0,\ell_i] \to [0,\ell'_i]$ for $i=1,2$ such that $\phi_1(\ell_1)=\ell'_1$ and $\phi_2(0)=0$ (we say that $\phi_1$ and $\phi_2$ are \textit{composable}). Define $\phi_1\ot\phi_2:[0,\ell_1+\ell_2]\to [0,\ell'_1+\ell'_2]$ as follows: 
	\[
	(\phi_1 \ot \phi_2)(t) = 
	\begin{cases}
		\phi_1(t) & \hbox{if } 0\leq t\leq \ell_1\\
		\phi_2(t-\ell_1)+\ell'_1 & \hbox{if } \ell_1\leq t\leq \ell_1+\ell_2.
	\end{cases}
	\] 
	We obtain a non-decreasing continuous map $\phi_1\ot\phi_2:[0,\ell_1+\ell_2]\to [0,\ell'_1+\ell'_2]$. When both $\phi_1$ and $\phi_2$ are surjective, i.e. when $\phi_1(0)=0$ and $\phi_2(\ell_2)=\ell'_2$, we recover the definition of Notation~\ref{nota:nondecreasing-surjective-ot}.
\end{nota}

\bp \label{lem-1-partial}
Let $\ell_1,\ell'_1,\ell_2,\ell'_2\geq 0$. Consider two composable non-decreasing continuous maps $\phi_i:[0,\ell_i] \to [0,\ell'_i]$ for $i=1,2$. Let $U$ be a topological space. Let $\gamma_i:[0,\ell'_i]\to U$ be two continuous maps with $i=1,2$. Then we have 
\[
\left(\gamma_1 *\gamma_2\right)(\phi_1\ot\phi_2) = (\gamma_1\phi_1)*(\gamma_2\phi_2).
\]
\ep

\bpf
For $0\leq t \leq \ell_1$, one has 
\begin{align*}
	\left(\gamma_1 *\gamma_2\right)(\phi_1\ot\phi_2)(t) &= \left(\gamma_1 *\gamma_2\right)(\phi_1(t)) \\
	& = \gamma_1(\phi_1(t)) \\
	& = \big((\gamma_1\phi_1)*(\gamma_2\phi_2)\big)(t),
\end{align*}
the first equality by definition of $\phi_1\ot \phi_2$, the second equality since $0\leq \phi_1(t)\leq \ell'_1$ and by definition of the Moore composition, and the third equality since $0\leq t \leq \ell_1$ and by definition of the Moore composition. 

For $\ell_1\leq t\leq \ell_1+\ell_2$, one has 
\begin{align*}
	\left(\gamma_1 *\gamma_2\right)(\phi_1\ot\phi_2)(t) &= \left(\gamma_1 *\gamma_2\right)(\phi_2(t-\ell_1)+\ell'_1) \\
	& = \gamma_2(\phi_2(t-\ell_1)) \\
	& = \big((\gamma_1\phi_1)*(\gamma_2\phi_2)\big)(t),
\end{align*}
the first equality by definition of $\phi_1\ot \phi_2$, the second equality since $\ell'_1\leq \phi_2(t-\ell_1)+\ell'_1\leq \ell'_1 + \ell'_2$ and by definition of the Moore composition, and the third equality since $\ell_1\leq t\leq \ell_1+\ell_2$ and by definition of the Moore composition.
\epf

\bp \label{lem0} \cite[Proposition~3.5]{Moore2}
Let $U$ be a topological space. Let $\gamma_i:[0,1]\to U$ be $n$ continuous maps with $1\leq i \leq n$ and $n\geq 1$. Let $\ell_i>0$ with $1\leq i \leq n$ nonzero real numbers with $\sum_i \ell_i = 1$. Then for all $\ell >0$, we have 
\[
\big((\gamma_1\mu_{\ell_1})*\dots *(\gamma_n\mu_{\ell_n})\big)\mu_{\ell} = (\gamma_1\mu_{\ell_1\ell})*\dots *(\gamma_n\mu_{\ell_n\ell}).
\]
\ep

\bd \cite[Definition~1.1]{reparam} \label{def_regular} Let $U$ be a topological space. The Moore path $\gamma:[0,L]\to U$ is \textit{regular} if $\gamma$ is constant or if for any in interval $[a,b]\subset [0,L]$ with $a\leq b$ such that the restriction $\gamma\rest_{[a,b]}$ is constant, one has $a=b$.
\ed 

The Moore composition of two non-constant regular paths in a Hausdorff space is therefore regular. The main property of regular paths for this paper is the following one: 

\bp \label{prop:regular_reparam}
Let $U$ be a Hausdorff topological space. Let $p,q:[0,1]\to U$ be two non-constant regular paths such that there exists $\phi\in \mathcal{M}(1,1)$ such that $p=q\phi$. Then $\phi$ is a homeomorphism from $[0,1]$ to itself. 
\ep

\bpf It is \cite[Proposition~3.8]{reparam}.
\epf

\bp \label{prop:local-inj-reg}
Let $U$ be a Hausdorff topological space. Let $\gamma:[0,L]\to U$ be a Moore path. Suppose that $\gamma$ is regular and that there exist $\eta_1,\eta_2\in \mathcal{M}(\ell,L)$ such that $\gamma\eta_1=\gamma\eta_2$. Then $\eta_1=\eta_2$. 
\ep

\bpf
It is \cite[Proposition~19]{Moore3} which is a Moore variant of \cite[Lemma~3.9]{reparam}.
\epf

Proposition~\ref{prop:regular_reparam} is used in Theorem~\ref{thm:T-homotopy}, and Proposition~\ref{prop:local-inj-reg} in Theorem~\ref{thm:decomposition_naturelle}. In both cases, the space $U$ is the underlying space of a cellular multipointed $d$-space. The latter is always Hausdorff: see Section~\ref{sec:directed-paths-cellular-case}. 

\section{From multipointed d-spaces to directed spaces}
\label{sec:from-multi-to-continuous}

\bd \cite[Definition~6]{Moore3} \label{def:multipointed-d-space} A \textit{multipointed $d$-space $X$} is a triple $(|X|,X^0,\P^{top}X)$ where
\begin{itemize}[leftmargin=*]
	\item The pair $(|X|,X^0)$ is a multipointed space. The space $|X|$ is called the \textit{underlying space} of $X$ and the set $X^0$ the \textit{set of states} of $X$.
	\item The set $\P^{top}X$ is a set of continuous maps from $[0,1]$ to $|X|$ called the \textit{execution paths}, satisfying the following axioms:
	\begin{itemize}
		\item For any execution path $\gamma$, one has $\gamma(0),\gamma(1)\in X^0$.
		\item Let $\gamma$ be an execution path of $X$. Then any composite $\gamma\phi$ with $\phi\in \mathcal{M}(1,1)$ is an execution path of $X$.
		\item Let $\gamma_1$ and $\gamma_2$ be two composable execution paths of $X$; then the normalized composition $\gamma_1 *_N \gamma_2$ is an execution path of $X$.
	\end{itemize}
\end{itemize}
A map $f:X\to Y$ of multipointed $d$-spaces is a map of multipointed spaces from $(|X|,X^0)$ to $(|Y|,Y^0)$ such that for any execution path $\gamma$ of $X$, the map $\P^{top}f:\gamma\mapsto f. \gamma$ is an execution path of $Y$. The category of multipointed $d$-spaces is denoted by $\ptop{\mathcal{M}}$. Let $\P_{\alpha,\beta}^{top} X = \{\gamma\in \P^{top}X\mid \gamma(0)=\alpha,\gamma(1)=\beta\}$.  
\ed

The category $\ptop{\mathcal{M}}$ is locally presentable, and in particular bicomplete, by e.g \cite[Proposition~6]{Moore3} which is an adaptation of \cite[Theorem~3.5]{mdtop}.

\bd \label{def:length-path}
Let $X$ be a multipointed $d$-space. The set of \textit{execution paths of length $\ell>0$} of $X$ from $\alpha\in X^0$ to $\beta\in X^0$ is the set $\P_{\alpha,\beta}^\ell X = \{\gamma\mu_\ell \mid \gamma\in \P_{\alpha,\beta}^{top}X\}$. It is equipped with the $\Delta$-kelleyfication of the relative topology with respect to the inclusion $\P_{\alpha,\beta}^\ell X \subset \ttop([0,\ell],|X|)$. Note that $\P_{\alpha,\beta}^1 X = \P^{top}_{\alpha,\beta} X$. 
\ed

\bp \label{prop:lengh-path}
Let $X$ be a multipointed $d$-space. Let $(\alpha,\beta,\gamma)\in X^0\p X^0\p X^0$. Let $\ell,\ell'>0$. Then there are the two following properties: 
\begin{itemize}
	\item For all $\gamma\in P_{\alpha,\beta}^\ell X$ and for all $\phi\in \mathcal{M}(\ell',\ell)$, one has $\gamma\phi\in \P_{\alpha,\beta}^{\ell'} X$.
	\item The Moore composition induces a continuous map 
	\[\P_{\alpha,\beta}^{\ell} X \p \P_{\beta,\gamma}^{\ell'} X \longrightarrow \P_{\alpha,\gamma}^{\ell+\ell'} X.\]
\end{itemize}
\ep

\bpf
The first assertion is \cite[Proposition~7]{Moore3}. The second assertion is \cite[Proposition~11]{Moore3}.
\epf

\begin{nota}
	Denote by $\mathcal{I}(\ell)$ the set of non-decreasing continuous maps from $[0,1]$ to $[0,\ell]$. Note that an element of $\mathcal{I}(\ell)$ can be a constant map.
\end{nota}

\bd \cite[Definition~1.1]{mg} \cite[Definition~4.1]{DAT_book} \label{def:directed_space}
A \textit{directed space} is a pair $X=(|X|,d(X))$ consisting of a topological space $|X|$ and a set $d(X)$ of continuous paths from $[0,1]$ to $|X|$ satisfying the following axioms:
\begin{itemize}
	\item $d(X)$ contains all constant paths;
	\item $d(X)$ is closed under normalized composition;
	\item $d(X)$ is closed under reparametrization by an element of $\mathcal{I}(1)$.
\end{itemize}
The space $|X|$ is called the \textit{underlying topological space} or the \textit{state space}. The elements of $d(X)$ are called \textit{directed paths}. A morphism of directed spaces is a continuous map between the underlying topological spaces which takes a directed path of the source to a directed path of the target. The category of directed spaces is denoted by $\ptop{}$. Write $\vec{P}(X)(u,v)$ for the space of directed paths of $X$ from $u$ to $v$ equipped with the $\Delta$-kelleyfication of the compact-open topology.
\ed

We use the terminology \textit{directed space} instead of the one of \textit{$d$-space} to avoid any confusion in this work with the \textit{multipointed $d$-spaces}. Likewise, directed spaces have directed paths and multipointed $d$-spaces have execution paths.

The category $\ptop{}$ is locally presentable. In particular, it is bicomplete. The proof is not strictly speaking in \cite{FR} although the latter paper contains the material. It suffices to start from a small relational universal strict Horn theory $\mathcal{T}$ axiomatizing $\top$ (without equality by \cite[Theorem~3.6]{FR} when $\top$ is the category of $\Delta$-generated spaces, and with equality by \cite[Proposition~B.18]{leftproperflow} when $\top$ is the category of $\Delta$-Hausdorff $\Delta$-generated spaces). The notion of directed space is axiomatized using the axioms expounded in the proof of \cite[Theorem~4.2]{FR}. The proof is complete thanks to \cite[Theorem~5.30]{TheBook}.

\bp \label{prop:cont}
Let $X$ be a multipointed $d$-space. Consider the set of paths $d(X)$ which consists of all constant paths and all Moore compositions of the form \[(\gamma_1\phi_1\mu_{\ell_1}) * \dots * (\gamma_n\phi_n\mu_{\ell_n})\] such that $\ell_1+\dots + \ell_n = 1$ where $\gamma_1,\dots,\gamma_n$ are execution paths of $X$ and $\phi_i\in \mathcal{I}(1)$ for $i=1,\dots,n$. Then the pair $(|X|,d(X))$ is a directed space denoted by $\cont(X)$. The mapping \[X\mapsto \cont(X)\] yields a well-defined functor from the category of multipointed $d$-spaces to the category of directed spaces.
\ep

\bpf By definition, the set $d(X)$ contains all constant paths. 

\paragraph{\textbf{(1)}} Let 
\[	 \Gamma = (\gamma_1\phi_1\mu_{\ell_1}) * \dots *(\gamma_n\phi_n\mu_{\ell_n}),\Gamma' =(\gamma'_1\phi'_1\mu_{\ell'_1}) * \dots *(\gamma'_{n'}\phi_{n'}\mu_{\ell_{n'}})\]
be two composable elements of $d(X)$. Then 
\begin{align*}
	\Gamma *_N \Gamma' & = (\Gamma\mu_{1/2}) * (\Gamma'\mu_{1/2}) \\
	& = \bigg((\gamma_1\phi_1\mu_{\ell_1/2}) * \dots (\gamma_n\phi_n\mu_{\ell_n/2}) * (\gamma'_1\phi'_1\mu_{\ell'_1/2}) * \dots (\gamma'_{n'}\phi_{n'}\mu_{\ell_{n'}/2})\bigg)
\end{align*}
the first equality by definition of the normalized composition, and the second equality by Proposition~\ref{lem0} applied to the continuous maps $\gamma_i\phi_i,\gamma'_i\phi'_i:[0,1]\to |X|$. Thus the normalized composition of two paths of $d(X)$ is a path of $d(X)$. 

\paragraph{\textbf{(2)}} Let $\phi\in \mathcal{I}(1)$. Let $\Gamma\in d(X)$. When $\phi(0)=\phi(1)$ or when $\Gamma$ is constant, $\Gamma\phi$ is constant as well and therefore $\Gamma\phi\in d(X)$. Assume now that $\Gamma$ is a non-constant path and that $0\leq\phi(0)<\phi(1) \leq 1$. Let $\Gamma = (\gamma_1\phi_1\mu_{\ell_1}) * \dots *(\gamma_n\phi_n\mu_{\ell_n})$. Using Proposition~\ref{lem-1} applied to the continuous maps $\gamma_i\phi_i:[0,1]\to |X|$, write 
\[
\Gamma = \bigg(\big(\gamma_1\phi_1\big) * \dots * \big(\gamma_n\phi_n\big)\bigg)\big(\mu_{\ell_1} \ot \dots \ot \mu_{\ell_n}\big).
\]
This implies 
\[
\Gamma\phi = \bigg(\big(\gamma_1\phi_1\big) * \dots * \big(\gamma_n\phi_n\big)\bigg)\underbrace{\big(\mu_{\ell_1} \ot \dots \ot \mu_{\ell_n}\big)\phi}_{\psi}.
\] 
One has $\psi\in \mathcal{I}(n)$. Since $\mu_{\ell_1} \ot \dots \ot \mu_{\ell_n}$ is a non-decreasing homeomorphism from $[0,1]$ to $[0,n]$, $\phi(0)<\phi(1)$ implies $\psi(0)< \psi(1)$. Let $r,s$ be the two unique integers such that 
\[
r= \max \{p\in \mathbb{N}\mid p\leq \psi(0)\} \leq \psi(0)< \psi(1) \leq s=\min \{p\in \mathbb{N}\mid \psi(1) \leq p\}.
\]
This implies that $r\leq \psi(0)<r+1$ and $s-1<\psi(1)\leq s$. From the equality $\psi(0)<\psi(1)$, we deduce that $r<s$, or equivalently $s -r \geq 1$. Then for all $t\in [0,1]$, one has by definition of the Moore composition
\[
\Gamma\phi(t) = \bigg(\big(\gamma_{r+1}\phi_{r+1}\big) * \dots * \big(\gamma_s\phi_s\big)\bigg){\big(\psi(t)-r\big)}.
\]
By definition of the Moore composition again, one obtains
\[
\Gamma\phi(t) = \begin{cases}
	\gamma_{r+1}\phi_{r+1}(\psi(t)-r) & \hbox{ if } \psi(0)\leq \psi(t)\leq r+1\\
	\gamma_{r+2}\phi_{r+2}(\psi(t)-(r+1)) & \hbox{ if } r+1\leq \psi(t)\leq r+2\\
	\dots\dots\dots\dots\dots\dots\dots \\
	\gamma_{s}\phi_{s}(\psi(t)-(s-1)) & \hbox{ if } s-1\leq \psi(t)\leq \psi(1).
\end{cases}
\]
Let $L_0=0$, $L_i\in [0,1]$ such that $\psi(L_i)=r+i$ for $1\leq i \leq s-r-1$ and $L_{s-r}=1$. The case $s-r=1$ means there are only two terms $L_0=0$ and $L_1=1$. In this case, $L_0<L_1$. In the case $s-r>1$, one has 
\begin{align*}
	& \psi(L_0) = \psi(0) < r+1 =\psi(L_1) \\
	& \psi(L_i) =r+i < r+i+1 = \psi(L_{i+1}) \hbox{ for }1\leq i<i+1\leq s-r-1\\
	& \psi(L_{s-r-1}) = s-1 <\psi(1)=\psi(L_{s-r}).
\end{align*}
In all cases, we deduce the strict inequalities \[\psi(L_0)<\psi(L_1)<\dots<\psi(L_{s-r}).\] This implies that $0=L_0<L_1<\dots<L_{s-r}=1$, $\psi$ being non-decreasing. We obtain
\[
\Gamma\phi(t) = \begin{cases}
	\gamma_{r+1}\phi_{r+1}(\psi(t)-r) & \hbox{ if } L_0\leq t\leq L_1\\
	\gamma_{r+2}\phi_{r+2}(\psi(t)-(r+1)) & \hbox{ if } L_1\leq t\leq L_2\\
	\dots\dots\dots\dots\dots\dots\dots \\
	\gamma_{s}\phi_{s}(\psi(t)-(s-1)) & \hbox{ if } L_{s-r-1}\leq t\leq L_{s-r}.
\end{cases}
\]
Let 
\begin{align*}
	&\overline{\phi}_{r+1}(u) = \phi_{r+1}(\psi(u+L_0)-r)  \hbox{ if } 0\leq u\leq L_1-L_0\\
	&\overline{\phi}_{r+2}(u) = \phi_{r+2}(\psi(u+L_1)-(r+1))  \hbox{ if } 0\leq u\leq L_2-L_1\\
	&\dots\dots\dots\dots\dots\dots\dots\dots\dots\dots\dots\dots\dots\dots\dots\dots\dots\dots \\
	&\overline{\phi}_{s}(u) = \phi_{s}(\psi(u+L_{s-r-1})-(s-1))  \hbox{ if } 0\leq u\leq L_{s-r}-L_{s-r-1}.
\end{align*}
Then, by definition of the Moore composition, we deduce 
\begin{align*}
	\Gamma\phi &= \big(\gamma_{r+1}\overline{\phi}_{r+1}\big) * \dots * \big(\gamma_{s}\overline{\phi}_{s}\big)\\
	&= \big(\gamma_{r+1}\overline{\phi}_{r+1}\mu^{-1}_{L_1-L_0}\mu_{L_1-L_0}\big)  * \dots * \big(\gamma_{s}\overline{\phi}_{s}\mu^{-1}_{L_{s-r}-L_{s-r-1}}\mu_{L_{s-r}-L_{s-r-1}}\big).
\end{align*}
By hypothesis, the maps $\phi_{r+1},\dots,\phi_s$ belong to $\mathcal{I}(1)$. Since $\psi$ is non-decreasing, belonging to $\mathcal{I}(n)$ by definition, the maps $\overline{\phi}_{r+1}\mu^{-1}_{L_1-L_0},\dots, \overline{\phi}_{s}\mu^{-1}_{L_{s-r}-L_{s-r-1}}$ belong to $\mathcal{I}(1)$ as well. We have proved that $\Gamma\phi\in d(X)$. 

\paragraph{\textbf{(3)}} Let $f:X\to Y$ be a map of multipointed $d$-spaces. Then $f$ takes constant paths of $X$ to constant paths of $X$ and an element of the form $(\gamma_1\phi_1\mu_{\ell_1}) * \dots *(\gamma_n\phi_n\mu_{\ell_n})$ to 
$(f\gamma_1\phi_1\mu_{\ell_1}) * \dots *(f\gamma_n\phi_n\mu_{\ell_n})$. Since $f\gamma_1,\dots,f\gamma_n$ are execution paths of $Y$, the mapping $X\mapsto \cont(X)$ induces a functor from multipointed $d$-spaces to directed spaces.
\epf

\bp\label{prop:Omega}
The mapping \[\discont:Y=(|Y|,d(Y)) \mapsto (|Y|,|Y|,d(Y))\] induces a full and faithful functor from the category of directed spaces $\ptop{}$ to the category of multipointed $d$-spaces $\ptop{\mathcal{M}}$. Moreover there is the equality \[\cont(\discont(X))=X\] for all directed spaces $X$.
\ep

\bpf
The set of execution paths of $\discont(Y)$ being closed under reparametrization by $\mathcal{M}(1,1)\subset \mathcal{I}(1)$ and closed under normalized composition, $\discont(Y)$ yields a well-defined multipointed $d$-space. 

From a map of directed spaces $f:Y\to Z$, we deduce a set map $f:Y^0=|Y| \to Z^0=|Z|$. Every execution path $\gamma$ of $\discont(Y)$ being a directed path of $d(Y)$, the map $f\gamma$ is a directed path of $d(Z)$, and therefore an execution path of $\discont(Z)$. Thus the mapping $Y\mapsto \discont(Y)$ induces a functor from the category of directed spaces to the one of multipointed $d$-spaces. 

Let $X$ be a directed space. By Proposition~\ref{prop:cont}, a directed path of $\cont(\discont(X))$ is a Moore composition of execution paths of $\discont(X)$ of various lengths such that the sum of these lengths is $1$. Thus, a directed path of $\cont(\discont(X))$ is an execution path of $\discont(X)$ by Proposition~\ref{prop:lengh-path}, and therefore a directed path of $X$. This implies that the identity of $|X|$ induces a map of directed spaces \[\cont(\discont(X)) \longrightarrow X.\] Since the identity is one-to-one, one obtains the inclusion $d(\cont(\discont(X))) \subset d(X)$. By construction of $\cont$, there is also the inclusion of sets $d(X)\subset d(\cont(\discont(X)))$. We obtain the equality \[\cont(\discont(X)) = X\] for all directed spaces $X$.

Let $f,g:Y\to Z$ be two maps of directed spaces such that $\discont(f)=\discont(g)$. Then $f=\cont(\discont(f))=\cont(\discont(g))=g$. Thus the set map \[\ptop{}(Y,Z) \to \ptop{\mathcal{M}}(\discont(Y),\discont(Z))\] is one-to-one. Let $g\in \ptop{\mathcal{M}}(\discont(Y),\discont(Z))$. Then $\cont(g):Y\to Z$ is a map of directed spaces. The maps of multipointed $d$-spaces $\discont(\cont(g)):\discont(Y)\to \discont(Z)$ and $g:\discont(Y)\to \discont(Z)$ induces the same continuous map $|g|:|Y|\to |Z|$. This implies that $\discont(\cont(g))=g$. Thus the set map \[\ptop{}(Y,Z) \to \ptop{\mathcal{M}}(\discont(Y),\discont(Z))\] is surjective and the functor $\discont:\ptop{}\to \ptop{\mathcal{M}}$ is full and faithful.
\epf

\bth \label{thm:cont-left-adjoint}
There is an adjunction $\cont\dashv \discont$. 
\eth

\bpf 
Let $X$ be a multipointed $d$-space. Let $Y$ be a directed space. By Proposition~\ref{prop:cont}, there is a unique set map $\ell$ from $\ptop{}(\cont(X),Y)$ to $\ptop{\mathcal{M}}(X,\discont(Y))$ such that the following diagram of sets 
\[
\begin{tikzcd}[row sep=3em, column sep=3em]
	\ptop{}(\cont(X),Y) \arrow[r,dashed,"\ell"] \arrow[d] & \ptop{\mathcal{M}}(X,\discont(Y)) \arrow[d] \\
	\top(|X|,|Y|) \arrow[r,equal] & \top(|X|,|Y|)
\end{tikzcd}
\]
is commutative where the vertical maps take a map to the corresponding map between the underlying spaces. Similarly, by Proposition~\ref{prop:cont}, there is a unique set map $\ell'$ from $\ptop{\mathcal{M}}(X,\discont(Y))$ to $\ptop{}(\cont(X),Y)$ such that the following diagram of sets 
\[
\begin{tikzcd}[row sep=3em, column sep=3em]
	\ptop{\mathcal{M}}(X,\discont(Y)) \arrow[r,dashed,"\ell'"] \arrow[d] &  \ptop{}(\cont(X),Y)\arrow[d] \\
	\top(|X|,|Y|) \arrow[r,equal] & \top(|X|,|Y|)
\end{tikzcd}
\]
is commutative. We deduce that $\ell\ell'$ is the unique map making the diagram
\[
\begin{tikzcd}[row sep=3em, column sep=3em]
	\ptop{\mathcal{M}}(X,\discont(Y)) \arrow[r,dashed,"\ell\ell'"]\arrow[d] &  \ptop{\mathcal{M}}(X,\discont(Y))\arrow[d] \\
	\top(|X|,|Y|) \arrow[r,equal] & \top(|X|,|Y|)
\end{tikzcd}
\]
commutative. Thus $\ell\ell'$ is the identity of $\ptop{\mathcal{M}}(X,\discont(Y))$. Similarly, we obtain that $\ell'\ell$ is the identity of $ \ptop{}(\cont(X),Y)$. Hence the proof is complete.
\epf

Since the categories $\ptop{}$ and $\ptop{\mathcal{M}}$ are locally presentable, it is possible to use a different argument to prove that the functor $\discont:\ptop{}\to \ptop{\mathcal{M}}$ is a right adjoint. From the natural bijection of sets \[\top([0,1],\limproj Z_i) \iso \limproj \top([0,1],Z_i)\] for all small diagrams of topological spaces $i\mapsto Z_i$, we see easily that the functor $\discont$ is limit-preserving. Since the set of directed paths (of execution paths resp.) of a colimit of directed spaces (of multipointed $d$-spaces resp.) is the set of finite compositions of directed paths (of execution paths resp.) of the components, the functor $\discont$ is colimit-preserving, and in particular accessible. Therefore by \cite[Theorem~1.66]{TheBook}, the functor $\discont$ is a right adjoint. 

On the other hand, we do not see any argument to prove without calculations that the set of paths $d(X)$ of Proposition~\ref{prop:cont} satisfies the axioms of directed paths. From the abstract argument above, we can only conclude that the left adjoint $\cont(X)$ contains as directed paths at least the set of paths $d(X)$ described in the statement of Proposition~\ref{prop:cont}. 

By the dual of the Special Adjoint Functor Theorem, the functor $\discont$ is also a left adjoint. The details are postponed to Appendix~\ref{sec:opcont} to avoid overloading the core of the text.

\bth \label{thm:reflection}
	The functor $\discont:\ptop{}\to \ptop{\mathcal{M}}$ induces an isomorphism between the category of directed spaces and a full reflective subcategory of the category of multipointed $d$-spaces. 
\eth

\bpf
By Proposition~\ref{prop:Omega}, there is the equality $\cont(\discont(X)) = X$. Finally, for every multipointed $d$-space of the form $\discont(X)$, one has $(\discont\cont)(\discont(X)) = \discont(X)$ by the previous equality. Hence the isomorphism of categories.
\epf

To justify the relevance of the functor $\cont:\ptop{\mathcal{M}}\to \ptop{}$ for DAT, we verify now that the directed spaces associated to the tame and non-tame realizations of a precubical set as a multipointed $d$-space are the same and that they coincide with the usual definition of the realization of a precubical set as a directed space. Proposition~\ref{prop:samerea} being not used in the sequel, it is possible to skip the reading of this part of the section.

We recall a few definitions: see e.g. \cite{DAT_book,MR4070250,RegularMoore}. Let $[n] = \{0,1\}^n$ for $n \geq 1$. Let $\{0,1\}^0=[0,1]^0=[0]=\{()\}$. Let $\delta_i^\alpha : [0,1]^{n-1} \rightarrow [0,1]^n$ be the continuous map defined for $1\leq i\leq n$ and $\alpha \in \{0,1\}$ by $\delta_i^\alpha(\epsilon_1, \dots, \epsilon_{n-1}) = (\epsilon_1,\dots, \epsilon_{i-1}, \alpha, \epsilon_i, \dots, \epsilon_{n-1})$. The small category $\square$ is the subcategory of the category of sets with the set of objects $\{[n],n\geq 0\}$ and generated by the \textit{coface maps} $\delta_i^\alpha$. A \textit{precubical set} $K$ is a presheaf over $\square$. An element $c$ of $K_n$ is called a \textit{$n$-cube} and we set $n=\dim(c)$. Let $\square[n]=\square(-,[n])$ for $n\geq 0$. The mapping $[n]\mapsto [0,1]^n$ yields a well-defined \textit{cocubical object} of $\top$, i.e. a functor from $\square$ to $\top$. The \textit{geometric realization} of the precubical set $K$ is the topological space 
\[
|K|_{geom} = \int^{[n]\in \square} K_n.[0,1]^n.
\]
A directed path of $[0,1]^n$ is a continuous map $\gamma:[0,1] \to [0,1]^n$ which is non-decreasing with respect to each axis of coordinates. It is \textit{tame} if $\gamma(0),\gamma(1)\in \{0,1\}^n$. A directed path of an $n$-cube $c$ of length $\ell \geq 0$ is a composite continuous map of the form $|c|_{geom}\gamma:[0,\ell] \to [0,1]^n\to |K|_{geom}$ such that $\gamma:[0,\ell]\to [0,1]^n$ is a directed path of length $\ell\geq 0$. A directed path in $|K|_{geom}$ of length $\ell\geq 0$ is a continuous path $[0,\ell] \to |K|_{geom}$ which is a Moore composition of the form $(|c_1|_{geom}\gamma_1) * \dots *(|c_p|_{geom}\gamma_p)$ with $p\geq 1$ such that $\ell=\ell_1+\dots+\ell_p$ where $\ell_i$ is the length of $\gamma_i$ for $1\leq i \leq p$. The choice of $c_1,\dots,c_p$ is not unique. This definition of a directed path in $|K|_{geom}$ makes sense anyway since the coface maps preserve the local ordering. The directed path $\gamma$ is \textit{tame} if each $\gamma_i$ for $1\leq i\leq p$ is tame (\cite[Section~2.9]{MR4070250}).

\begin{figure}
	\begin{tikzpicture}
		\draw (0,0) -- (0,2) -- (2,2) -- (2,0) -- (0,0);
		\draw (2,2) -- (4,2) -- (4,0) -- (2,0);
		\draw (2,2) -- (2,4) -- (4,4) -- (4,2);
		\draw (5,0) -- (5,2) -- (7,2) -- (7,0) -- (5,0);
		\draw (7,2) -- (9,2) -- (9,0) -- (7,0);
		\draw (7,2) -- (7,4) -- (9,4) -- (9,2);	
		\draw[->, line width=0.4mm, color=dark-red, smooth] plot coordinates {(0,0) (0.5,1.5) (2,2)};
		\draw[->, line width=0.4mm, color=dark-red, smooth] plot coordinates {(2,2) (3.5,2.5) (4,4)};		
		\draw[->, line width=0.4mm, color=dark-red, smooth] plot coordinates {(5,0) (6,0.7) (7,1)};
		\draw[->, line width=0.4mm, color=dark-red, smooth] plot coordinates {(7,1) (8,1.3) (8.2,2)};
		\draw[->, line width=0.4mm, color=dark-red, smooth] plot coordinates {(8.2,2) (8.5,3) (9,4)};	
	\end{tikzpicture}
	\caption{The left directed path is tame, the right directed path is not tame.}
	\label{tame-nontame}
\end{figure}
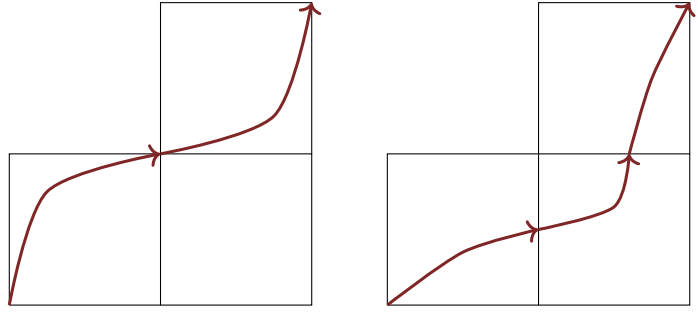 

\bp \label{cocubical-reg-cube} 
Let $n\geq 1$. The following data assemble into a multipointed $d$-space denoted by $|\square[n]|^t$:
\begin{itemize}
	\item The underlying space is the topological $n$-cube $[0,1]^n$.
	\item The set of states is $\{0,1\}^n\subset [0,1]^n$.
	\item The set of execution paths from $\underline{a}$ to $\underline{b}$ with $\underline{a}<\underline{b} \in \{0<1\}^n$ is the set of directed paths $[0,1]\to [0,1]^n$ from $\underline{a}$ to $\underline{b}$.
	\item The set of execution paths from $\underline{a}$ to $\underline{b}$ when $\underline{a}\geq \underline{b}$ or $\{\underline{a},\underline{b}\}$ is an incomparable pair is empty.
\end{itemize}
Let $|\square[0]|^t = \{()\}$. The mapping $[n]\mapsto |\square[n]|^t$ yields a well-defined cocubical objects of $\ptop{\mathcal{M}}$. 
\ep

\bpf
It is an adaptation of \cite[Proposition~2.22]{RegularMoore} in the non-regular case.
\epf

This leads to the definition (it is an adaptation to the non-regular case of \cite[Definition~2.23]{RegularMoore}): 

\bd \label{t_rea}
Let $K$ be a precubical set. The \textit{tame realization} of $K$ is the multipointed $d$-space 
\[
|K|^t = \int^{[n]\in \square} K_n.|\square[n]|^t.
\]
Note that the underlying space of $|K|^t$ is the topological space $|K|_{geom}$.
\ed

This yields a colimit-preserving functor from precubical sets to multipointed $d$-spaces by \cite[Proposition~2.15]{RegularMoore}. The chosen name \textit{tame realization} is because all execution paths of $|K|^t$ are tame. Another realization functor from precubical sets to multipointed $d$-spaces can be defined as follows (it is an adaptation to the non-regular case of \cite[Definition~7.1]{RegularMoore}): 

\bd \label{rea}
Let $K$ be a precubical set. The \textit{realization} $|K|$ of $K$ is the multipointed $d$-space having  as underlying space $|K|_{geom}$, the set of states $K_0$, and such that the set of execution paths from $\alpha$ to $\beta$ consists of the \textit{non-constant} directed paths from $\alpha$ to $\beta$ in $|K|_{geom}$.
\ed

We obtain the important fact: 

\bp \label{prop:samerea}
Let $K$ be a precubical set. The identity of $|K|_{geom}$ induces a map of multipointed $d$-spaces $|K|^t\longrightarrow |K|$ from the tame realization to the realization of $K$. This map induces an isomorphism of directed spaces 
\[
\cont(|K|^t) \iso \cont(|K|).
\]
\ep

\bpf
The existence of the map of multipointed $d$-spaces $|K|^t\longrightarrow |K|$ is due to the fact that the identity of $|K|_{geom}$ induces a one-to-one set map from the set of non-constant tame directed paths to the set of non-constant directed paths between two vertices in the geometric realization of $K$. Thus, by Proposition~\ref{prop:cont}, and since the underlying map is one-to-one, there is a one-to-one set map from $d(|K|^t)$ to $d(|K|)$. To prove the isomorphism $\cont(|K|^t) \iso \cont(|K|)$, it remains to prove that every directed path of $\cont(|K|)$ is a directed path of $\cont(|K|^t)$. By Proposition~\ref{prop:cont}, it suffices to prove that every execution path of $|K|$ is a directed path of $\cont(|K|^t)$. Every execution path $\gamma$ of $|K|$ is a Moore composition of the form $\gamma = (|c_1|_{geom}\gamma_1) * \dots *(|c_p|_{geom}\gamma_p)$ such that each directed path $\gamma_i$ in $[0,1]^{\dim(c_i)}$ is of length $\ell_i>0$ with $\ell_1+\dots+\ell_p=1$. Each directed path $\gamma_i$ can be precomposed with a directed path from $(0,\dots,0)$ to $\gamma_i(0)$ in $[0,1]^{\dim(c_i)}$ and postcomposed with a directed path from $\gamma_i(\ell_i)$ to $(1,\dots,1)$ in $[0,1]^{\dim(c_i)}$ to obtain a directed path $\overline{\gamma_i}$ of length $L_i\geq \ell_i>0$ from $(0,\dots,0)$ to $(1,\dots,1)$ in $[0,1]^{\dim(c_i)}$. There exists a non-decreasing map $\phi_i:[0,\ell_i]\to [0,L_i]$ such that $\gamma_i=\overline{\gamma_i}\phi_i$ for each $i\in \{1,\dots,p\}$. We obtain \[\gamma = \big(|c_1|_{geom}\overline{\gamma_1}\mu^{-1}_{L_1}\mu_{L_1}\phi_1\mu^{-1}_{\ell_1}\mu_{\ell_1}\big) * \dots *\big(|c_p|_{geom}\overline{\gamma_p}\mu^{-1}_{L_p}\mu_{L_p}\phi_p\mu^{-1}_{\ell_p}\mu_{\ell_p}\big).\] Each $\overline{\gamma_i}\mu^{-1}_{L_i}$ being tame by definition, each $|c_i|_{geom}\overline{\gamma_i}\mu^{-1}_{L_i}$ is an execution path of $|K|^t$. Each $\mu_{L_i}\phi_i\mu^{-1}_{\ell_i}$ for $1\leq i \leq p$ belongs to $\mathcal{I}(1)$. Using Proposition~\ref{prop:cont}, we deduce that $\gamma$ is a directed path of $\cont(|K|^t)$. 
\epf

As a consequence of Proposition~\ref{prop:cont} and Proposition~\ref{prop:samerea}, the directed space $\cont(|K|^t) \iso \cont(|K|)$ is nothing else but the usual directed space associated with a precubical set as defined e.g. in \cite{DAT_book}.

\section{The directed paths of a cellular multipointed d-space}
\label{sec:directed-paths-cellular-case}

The \textit{topological globe of $Z$}, which is denoted by $\globM(Z)$, is the multipointed $d$-space defined as follows
\begin{itemize}
	\item the underlying topological space is the quotient space \[\frac{\{{0},{1}\}\sqcup (Z\p[0,1])}{(z,0)=(z',0)={0},(z,1)=(z',1)={1}}\]
	\item the set of states is $\{{0},{1}\}$
	\item the set of execution paths is the set of continuous maps \[\{\delta_z\phi\mid \phi\in \mathcal{M}(1,1),z\in  Z\}\]
	with $\delta_z(t) = (z,t)$.	It is equal to the underlying set of the space $Z \p \mathcal{M}(1,1)$.
\end{itemize}
In particular, $\globM(\varnothing)$ is the multipointed $d$-space $\{{0},{1}\} = (\{{0},{1}\},\{{0},{1}\},\varnothing)$. The \textit{directed segment} is the multipointed $d$-space $\vI^{top}=\globM(\{0\})$. 

Let $n\geq 1$. Denote by $\mathbf{D}^n = \{(x_1,\dots,x_n)\in \mathbb{R}^n, x_1^2+\dots + x_n^2 \leq 1\}$ the $n$-dimensional disk, and by $\mathbf{S}^{n-1} = \{(x_1,\dots,x_n)\in \mathbb{R}^n, x_1^2+\dots + x_n^2 = 1\}$ the $(n-1)$-dimensional sphere. By convention, let $\mathbf{D}^{0}=\{0\}$ and $\mathbf{S}^{-1}=\varnothing$.

A \textit{cellular multipointed $d$-space} is a multipointed $d$-space which is the transfinite composition of a \textit{cellular decomposition}. The latter consists of a colimit-preserving functor \[X:\lambda \longrightarrow \ptop{\mathcal{M}}\] where $\lambda$ is a transfinite ordinal viewed as a small category thanks to its poset structure such that
\begin{itemize}
	\item The multipointed $d$-space $X_0$ is a set, in other terms $X_0=(X^0,X^0,\varnothing)$ for some set $X^0$.
	\item For all $\nu<\lambda$, there is a pushout diagram of multipointed $d$-spaces 
	\[
	\begin{tikzcd}[row sep=3em, column sep=3em]
		\globM(\mathbf{S}^{n_\nu-1}) \arrow[d] \arrow[r,"g_\nu"] & X_\nu \arrow[d] \\
		\globM(\mathbf{D}^{n_\nu}) \arrow[r,"\widehat{g_\nu}"] &  \cocartesian X_{\nu+1}
	\end{tikzcd}
	\]
	with $n_\nu \geq 0$. 
\end{itemize}
The cellular multipointed $d$-spaces are the cellular objects of the combinatorial model category whose definition is recalled in Section~\ref{sec:counter-example}. Note that the cofibration $R:\{0,1\}\to \{0\}$ is not required to characterize the class of cellular objects. 

For the sequel, for any cellular multipointed $d$-space, a cellular decomposition is implicitly fixed. A cellular multipointed $d$-space can be the transfinite composition of different cellular decompositions. Indeed, for example two cellular decompositions of $\mathbf{D}^2$ like 
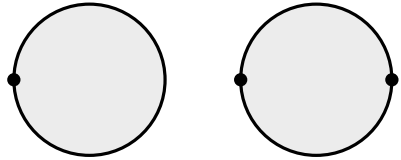
\[
\begin{tikzpicture}[black,pn/.style={circle,inner sep=0pt,minimum width=5pt,fill=black}]
	\fill[color=gray!15] (0,0) circle (1cm);
	\fill[color=gray!15] (3,0) circle (1cm);
	\draw (2,0) node[pn] {};
	\draw (4,0) node[pn] {};
	\draw (-1,0) node[pn] {};
	\draw [very thick] (0,0) circle (1cm);
	\draw [very thick] (3,0) circle (1cm);
\end{tikzpicture}
\]
will give rise to two cellular decompositions of $\globM(\mathbf{D}^2)$.

\begin{nota}
	Let $X_\lambda = \liminj_{\nu<\lambda} X_\nu$ be the transfinite composition.
\end{nota}

\begin{figure}
	\begin{tikzpicture}[black,scale=2.5,pn/.style={circle,inner sep=0pt,minimum width=4pt,fill=dark-red}]
		\draw (0,0) node[pn] {} node[black,below left] {$(0\,,0)$};
		\draw (1,0) node[pn] {} node[black,below right] {$(1\,,0)$};
		\draw (0,1) node[pn] {} node[black,above left] {$(0\,,1)$};
		\draw (1,1) node[pn] {} node[black,above right] {$(1\,,1)$};
		\draw[->] [very thick] (0.2,0.2) -- (0.5,0.5);
		\draw[->] [very thick] (0.5,0.5) -- (0.8,0.2);
		\draw[->] (0,0) -- (1,1);
		\draw[->] (0,1) -- (1,0); 
	\end{tikzpicture}
	\caption{The fake crossing}
	\label{counterexample}
\end{figure}
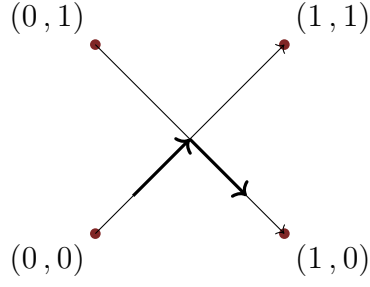

Consider the non-cellular multipointed $d$-space $X$ depicted in Figure~\ref{counterexample}. It is defined as follows. The underlying space $|X|$ is the topological space  $|X| = \{(u,u)\mid u\in [0,1]\} \cup \{(u,1-u)\mid u\in [0,1]\}$. Let $X^0 = \{(0,0),(0,1),(1,0),(1,1)\}$. Let $\P^{top}_{(0,0),(1,1)}X = \{t\mapsto (\phi(t),\phi(t))\mid \phi\in \mathcal{M}(1,1)\}$, $\P^{top}_{(0,1),(1,0)}X = \{t\mapsto (\phi(t),1-\phi(t))\mid \phi\in \mathcal{M}(1,1)\}$ and $\P^{top}_{\alpha,\beta}X=\varnothing$ otherwise. Then the composite of the two thick directed paths of Figure~\ref{counterexample} is not of the form $\gamma\phi$ where $\gamma$ is an execution path of $X$ and where $\phi\in \mathcal{I}(1)$. We prove in Theorem~\ref{continuous} that all directed paths of a cellular multipointed $d$-space are of this form.

\begin{nota}
	In the whole section, $X_\lambda$ stands for a cellular multipointed $d$-space like above.
\end{nota}

The underlying topological space $|X_\lambda|$ is Hausdorff by \cite[Proposition~15]{Moore3}. For all $\nu\leq \lambda$, there is the equality $X_\nu^0=X^0$. Denote by \[c_\nu = |\globM(\mathbf{D}^{n_\nu})|\backslash |\globM(\mathbf{S}^{n_\nu-1})|\] the $\nu$-th cell of $X_\lambda$. It is called a \textit{globular cell}. Like in the usual setting of CW-complexes, $\widehat{g_\nu}$ induces a homeomorphism from $c_\nu$ to $\widehat{g_\nu}(c_\nu)$ equipped with the relative topology. The map $\widehat{g_{\nu}}: \globM(\mathbf{D}^{n_\nu})\to X_\lambda$ is called the \textit{attaching map} of the globular cell $c_\nu$. The state $\widehat{g_\nu}(0)\in X^0$ ($\widehat{g_\nu}(1)\in X^0$ resp.)  is called the \textit{initial (final resp.) state} of $c_\nu$ and is denoted by $c_\nu^-$ ($c_\nu^+$ resp.). The integer $n_\nu+1$ is called the \textit{dimension} of the globular cell $c_\nu$. It is denoted by $\dim c_\nu$. The states of $X^0$ are also called the \textit{globular cells of dimension $0$}. By convention, a state of $X^0$ viewed as a globular cell of dimension $0$ is equal to its initial state and to its final state. Thus, for $\alpha\in X^0$, one has $\alpha=\alpha^+=\alpha^-$. The set of globular cells of $X_\lambda$ is denoted by $\mathcal{C}(X_\lambda)$. The set of globular cells of dimension $n\geq 0$ of $X_\lambda$ is denoted by $\mathcal{C}_n(X_\lambda)$. In particular, $\mathcal{C}_0(X_\lambda)=X^0$.

\begin{nota}
	Let $c \neq d$ be two globular cells of $X_\lambda$. The notation $c\preceq d$ means that either $c=d^-$ (which implies $\dim(d)\geq 1$) or $c^+=d$ (which implies $\dim(c)\geq 1$).
\end{nota}

\bd \label{def:dt}
A \textit{discrete trace} is a finite sequence of globular cells $[c_1,\dots,c_n]$ of $\mathcal{C}(X_\lambda)$ such that $c_1\preceq \dots \preceq c_n$. 
\ed

\bp \label{prop:dt-constant}
Let $x\in |X_\lambda|$. There exists a unique globular cell $\dt(x)$ of $X_\lambda$ such that $x\in \dt(x)$.
\ep

\bpf
The underlying topological space $|X_\lambda|$ satisfies as a set the equality of sets
\[
|X_\lambda| = \coprod_{c\in \mathcal{C}(X_\lambda)} c.
\]
\epf

\bp \label{prop:example_reg} \cite[Proposition~18]{Moore3}
Consider a globular cell $c_\nu$ of $X_\lambda$. Let $z\in \mathbf{D}^{n_\nu}\backslash \mathbf{S}^{n_\nu-1}$. The execution path $\widehat{g_{\nu}}\delta_{z}$ is regular.
\ep

Note that $\widehat{g_{\nu}}\delta_{z}$ is not necessarily regular if $z\in \mathbf{S}^{n_\nu-1}$.

\bth \label{thm:decomposition_naturelle} 
Let $\gamma$ be an execution path of $X_\lambda$. There is a unique decomposition of $\gamma$ of the form $\natgl(\gamma)\phi$ with $\natgl(\gamma) = (\widehat{g_{\nu_1}}\delta_{z_1})*\dots * (\widehat{g_{\nu_n}}\delta_{z_n})$, $n\geq 1$, $\nu_i<\lambda$ and $z_i\in \mathbf{D}^{n_{\nu_i}}\backslash \mathbf{S}^{n_{\nu_i}-1}$ for $1\leq i\leq n$, and $\phi\in \mathcal{M}(1,n)$. Finally, the path $\natgl(\gamma)$ is a regular execution path of length $n$.
\eth

\bpf
The existence of the decomposition $\natgl(\gamma)\phi$ is a consequence of \cite[Theorem~6]{Moore3}. Consider a second decomposition of the form $\natgl'(\gamma)\phi'$ of $\gamma$. By \cite[Theorem~6]{Moore3}, one has $\natgl(\gamma)=\natgl'(\gamma)$. We deduce that $\gamma=\natgl(\gamma)\phi=\natgl(\gamma)\phi'$. Since a Moore composition of regular paths in the Hausdorff space $|X_\lambda|$ is regular, the path $\natgl(\gamma)$ is regular by Proposition~\ref{prop:example_reg}. Thus by Proposition~\ref{prop:local-inj-reg}, one obtains $\phi=\phi'$. Finally, the path $\natgl(\gamma)$ is an execution path of length $n$ by Proposition~\ref{prop:lengh-path}. 
\epf

\bd \cite[Definition~11]{Moore3}
The execution path $\natgl(\gamma)$ of Theorem~\ref{thm:decomposition_naturelle} is called the \textit{globular naturalization} of $\gamma$. The sequence $\carrier(\gamma) = [c_{\nu_1},\dots,c_{\nu_n}]$ of Theorem~\ref{thm:decomposition_naturelle} is called the \textit{carrier} of $\gamma$. The integer $n$ is called the \textit{length} of the carrier. 
\ed

\bth \label{continuous} 
A continuous map from $[0,1]$ to $|X_\lambda|$ is a directed path of $X_\lambda$, i.e. of $\cont(X_\lambda)$, if and only if it is of the form $\gamma\phi$ where $\gamma$ is an execution path of $X_\lambda$ or a constant path and where $\phi\in \mathcal{I}(1)$. 
\eth

\bpf
It suffices to prove that for any execution path $\gamma$ and $\gamma'$ of $X_\lambda$ and any $\phi,\phi'\in \mathcal{I}(1)$ such that $\gamma\phi$ and $\gamma'\phi'$ are two composable paths of $|X_\lambda|$, i.e. $\gamma\phi(1)= \gamma'\phi'(0)$, the normalized composition $(\gamma\phi) *_N (\gamma'\phi')$ is of the form $\gamma''\phi''$ where $\gamma''$ is an execution path of $X_\lambda$ and $\phi''\in \mathcal{I}(1)$. Using Theorem~\ref{thm:decomposition_naturelle}, write 
\begin{align*}
	& \carrier(\gamma)=[c_{\nu_1},\dots,c_{\nu_n}], \\
	& \carrier(\gamma')=[c_{\nu'_1},\dots,c_{\nu'_{n'}}], \\
	& \natgl(\gamma) = (\widehat{g_{\nu_1}}\delta_{z_1})*\dots * (\widehat{g_{\nu_n}}\delta_{z_n}),\\
	& \natgl(\gamma') = (\widehat{g_{\nu'_1}}\delta_{z'_1})*\dots * (\widehat{g_{\nu'_{n'}}}\delta_{z'_{n'}}),\\
	& \gamma=\natgl(\gamma)\psi\hbox{ with }\psi\in \mathcal{M}(1,n),\\
	& \gamma'=\natgl(\gamma')\psi'\hbox{ with }\psi'\in \mathcal{M}(1,n').
\end{align*}
There are two mutually exclusive cases: (1) $\gamma\phi(1)= \gamma'\phi'(0) \in X^0$, (2) $\gamma\phi(1)= \gamma'\phi'(0) \notin X^0$.

\paragraph{\textbf{(1)}} One has $\gamma\phi(1)= \natgl(\gamma)(\psi\phi(1)) \in X^0$. We deduce that $\psi\phi(1)\in \{0,1,\dots,n\}$, the execution path $\natgl(\gamma)$ being regular by Theorem~\ref{thm:decomposition_naturelle}. Similarly, one has $\gamma'\phi'(0)= \natgl(\gamma')(\psi'\phi'(0)) \in X^0$. We deduce that $\psi'\phi'(0)\in \{0,1,\dots,n'\}$. We can suppose without lack of generality that $\psi\phi(1)=n$ and $\psi'\phi'(0) = 0$ by shortening $\gamma$ and $\gamma'$ (which implies in particular to remove some globular cells from their carrier). This implies that the non-decreasing continuous maps $\psi\phi\mu_{1/2}:[0,1/2]\to [0,n]$ and $\psi'\phi'\mu_{1/2}:[0,1/2]\to [0,n']$ are composable in the sense of Notation~\ref{nota:nondecreasing-continuous-ot}. Then one obtains the sequence of equalities
\begin{align*}
	(\gamma\phi) *_N (\gamma'\phi') &= (\gamma\phi\mu_{1/2}) * (\gamma'\phi'\mu_{1/2}) \\
	&= \big(\natgl(\gamma)\psi\phi\mu_{1/2}\big) * \big(\natgl(\gamma')\psi'\phi'\mu_{1/2}\big) \\
	&= \big(\natgl(\gamma)*\natgl(\gamma')\big)\big((\psi\phi\mu_{1/2}) \ot (\psi'\phi'\mu_{1/2})\big),
\end{align*}
the first one by definition of $*_N$, the second one by trivial substitution, and the third one by Proposition~\ref{lem-1-partial} and since $\psi\phi\mu_{1/2}$ and $\psi'\phi'\mu_{1/2}$ are composable. This implies that \[(\gamma\phi) *_N (\gamma'\phi') = \gamma''\phi''\] where $\gamma''$ is an execution path of $X_\lambda$ by Theorem~\ref{thm:decomposition_naturelle} and Proposition~\ref{prop:lengh-path} and $\phi''\in \mathcal{I}(1)$ with 
\[
\gamma'' = \big(\natgl(\gamma)*\natgl(\gamma')\big) \mu_{m+n}^{-1}
\]
and 
\[
\phi'' = \mu_{m+n}\bigg((\psi\phi\mu_{1/2}) \ot (\psi'\phi'\mu_{1/2})\bigg).
\]

\paragraph{\textbf{(2)}} In this case, $\psi\phi(1)\in [0,n]\backslash \{0,1,\dots,n\}$ and $\psi'\phi'(0)\in [0,n']\backslash \{0,1,\dots,n'\}$, the execution paths $\natgl(\gamma)$ and $\natgl(\gamma')$ being regular by Theorem~\ref{thm:decomposition_naturelle}. We can suppose without lack of generality that $n-1<\psi\phi(1)<n$ and $0<\psi'\phi'(0) <1$ by shortening $\gamma$ and $\gamma'$ (which implies in particular to remove some globular cells from their carrier). From the equality $\gamma\phi(1)= \gamma'\phi'(0) \notin X^0$, we deduce that $\gamma\phi(1)= \gamma'\phi'(0) \in c_{\nu_n} \cap c_{\nu'_1}$. From the set bijection \[|X_\lambda| = X^0 \sqcup \coprod_{\nu<\lambda} c_\nu,\] 
we obtain that $c_{\nu_n} = c_{\nu'_1}$ and $\nu_n=\nu'_1$. From the equality $\gamma\phi(1)= \gamma'\phi'(0)$, we then deduce that $z_n=z'_1$. We obtain the equality $\widehat{g_{\nu_n}}\delta_{z_n} = \widehat{g_{\nu'_1}}\delta_{z'_1}$. By definition of the Moore composition, we obtain the equality 
\[
\widehat{g_{\nu_n}}(z_n,\psi\phi(1)-(n-1)) = \widehat{g_{\nu'_1}}(z'_1,\psi'\phi'(0)).
\]
From the fact that $\widehat{g_{\nu_n}}\delta_{z_n} = \widehat{g_{\nu'_1}}\delta_{z'_1}$ is regular by Proposition~\ref{prop:example_reg}, we also deduce the equality 
\[
\psi\phi(1)-(n-1)=\psi'\phi'(0)\in ]0,1[.
\]
Thanks to the equalities $\widehat{g_{\nu_n}}\delta_{z_n} = \widehat{g_{\nu'_1}}\delta_{z'_1}$ and $\psi\phi(1)-(n-1)=\psi'\phi'(0)$, we introduce two new objects, a Moore path $\Gamma$ of $|X_\lambda|$ of length $n+n'-1$  defined in two equivalent ways by
\begin{multline*}
	\Gamma = (\widehat{g_{\nu_1}}\delta_{z_1})*\dots * (\widehat{g_{\nu_n}}\delta_{z_n}) * (\widehat{g_{\nu'_2}}\delta_{z'_2})*\dots * (\widehat{g_{\nu'_{n'}}}\delta_{z'_{n'}})\\ = (\widehat{g_{\nu_1}}\delta_{z_1})*\dots * (\widehat{g_{\nu_{n-1}}}\delta_{z_{n-1}}) * (\widehat{g_{\nu'_1}}\delta_{z'_1})*\dots * (\widehat{g_{\nu'_{n'}}}\delta_{z'_{n'}})
\end{multline*}
and a non-decreasing continuous map $\Phi:[0,1]\to [0,n+n'-1]$ defined by 
\[
\Phi(t) = 
\begin{cases}
	\psi(\phi(2t)) & \hbox{ if } 0\leq t\leq \frac{1}{2}\\
	\psi'(\phi'(2t-1)) + (n-1) & \hbox{ if } \frac{1}{2}\leq t\leq 1
\end{cases}
\]
The continuous path $\Gamma$ is an execution path of $X_\lambda$ of length $n+n'-1$ by Proposition~\ref{prop:lengh-path}. One obtains 
\[
\Gamma\Phi(t) =\begin{cases}
	\widehat{g_{\nu_p}}(z_p,\psi\phi(2t)-(p-1)) & \hbox{ if }t\leq \frac{1}{2}\hbox{ and }p-1\leq \psi\phi(2t)\leq p\\
	\widehat{g_{\nu'_p}}(z'_p,\psi'\phi'(2t-1)-(p-1)) & \hbox{ if }t\geq \frac{1}{2}\hbox{ and }p-1\leq \psi'\phi'(2t-1)\leq p.
\end{cases}
\]
the first equality since $\Gamma = (\widehat{g_{\nu_1}}\delta_{z_1})*\dots * (\widehat{g_{\nu_n}}\delta_{z_n}) * (\widehat{g_{\nu'_2}}\delta_{z'_2})*\dots * (\widehat{g_{\nu'_{n'}}}\delta_{z'_{n'}})$ and the second equality since $\Gamma = (\widehat{g_{\nu_1}}\delta_{z_1})*\dots * (\widehat{g_{\nu_{n-1}}}\delta_{z_{n-1}}) * (\widehat{g_{\nu'_1}}\delta_{z'_1})*\dots * (\widehat{g_{\nu'_{n'}}}\delta_{z'_{n'}})$. Let 
\[\Gamma':=(\gamma\phi) *_N (\gamma'\phi') = (\natgl(\gamma)\psi\phi\mu_{1/2}) * (\natgl(\gamma')\psi'\phi'\mu_{1/2}).\] 
One obtains by definition of the Moore composition
\[
\Gamma'(t) =\begin{cases}
	\widehat{g_{\nu_p}}(z_p,\psi\phi(2t)-(p-1)) & \hbox{ if }t\leq \frac{1}{2}\hbox{ and }p-1\leq \psi\phi(2t)\leq p\\
	\widehat{g_{\nu'_p}}(z'_p,\psi'\phi'(2t-1)-(p-1)) & \hbox{ if }t\geq \frac{1}{2}\hbox{ and }p-1\leq \psi'\phi'(2t-1)\leq p.
\end{cases}
\]
We deduce that \[(\gamma\phi) *_N (\gamma'\phi') = \Gamma\Phi.\] This implies that \[(\gamma\phi) *_N (\gamma'\phi') = \gamma''\phi''\] where $\gamma''$ is an execution path of $X_\lambda$ by Proposition~\ref{prop:lengh-path} and $\phi''\in \mathcal{I}(1)$ with \[\gamma'' = \Gamma \mu_{n+n'-1}^{-1}\] and \[\phi'' = \mu_{n+n'-1}\Phi.\]
\epf

As a first application of Theorem~\ref{continuous}, we can define the discrete trace of a directed path of a cellular multipointed $d$-spaces as follows. The formulation of Proposition~\ref{constr-dt} is close to \cite[Definition~2.20]{zbMATH02231448}. Another application of Theorem~\ref{continuous} is given in Theorem~\ref{thm:T-homotopy}.

\bth \label{constr-dt}
Let $\gamma$ be a directed path of $X_\lambda$. There exists a unique discrete trace $\dt(\gamma)=[c_1,\dots,c_n]$ and a sequence of real numbers $0=t_0\leq t_1\leq \dots \leq t_n=1$ such that 
\begin{itemize}
	\item $c_i\neq c_{i+1}$,
	\item $t\in [t_{i-1},t_{i}] \Rightarrow \gamma(t)\in \widehat{c_i}$,
	\item $t\in ]t_{i-1},t_{i}[ \Rightarrow \dt(\gamma(t)) = c_i$ ,
	\item $\dt(\gamma(t_i)) \in \{c_i,c_{i+1}\}$,
	\item $\dt(\gamma(0)) = c_1$ and $\dt(\gamma(1))=c_n$.
\end{itemize} 
\eth

\bpf
Suppose at first that the directed path $\gamma$ is an execution path of $X_\lambda$. Using Theorem~\ref{thm:decomposition_naturelle}, let $\carrier(\gamma)=[c_{\nu_1},\dots,c_{\nu_n}]$, $\natgl(\gamma) = (\widehat{g_{\nu_1}}\delta_{z_1})*\dots * (\widehat{g_{\nu_n}}\delta_{z_n})$ and $\gamma=\natgl(\gamma)\psi$ with $\psi\in \mathcal{M}(1,n)$. Then one has \[\dt(\gamma) = [\widehat{g_{\nu_1}}(\bullet,0),c_{\nu_1},\widehat{g_{\nu_2}}(\bullet,0),c_{\nu_2},\dots,\widehat{g_{\nu_n}}(\bullet,0),c_{\nu_1},\widehat{g_{\nu_n}}(\bullet,1)],\] the symbol $\bullet$ meaning here that its value does not matter. Note that $\carrier(\gamma) \neq \dt(\gamma)$, the discrete trace of $\gamma$ containing also $0$-dimensional globular cells. In the general case, a directed path of $X_\lambda$ is either of the form $\gamma\phi$ with $\phi\in \mathcal{I}(1)$ or a constant path by Theorem~\ref{continuous}. In the first case, $\dt(\gamma\phi)$ will be a subsequence of $\dt(\gamma)$. In the case of a constant path $x$, the discrete trace is given by $[\dt(x)]$ of Proposition~\ref{prop:dt-constant}. 
\epf

\bd
Let $\gamma$ be a directed path of $X_\lambda$. The finite sequence of globular cells $\dt(\gamma) = [c_1,\dots,c_n]$ of Theorem~\ref{constr-dt} is called the \textit{discrete trace} of $\gamma$. Note the abuse of notation, $\dt(x)$ meaning either the unique globular cell containing $x\in |X_\lambda|$ and the discrete trace of the constant path $x$.
\ed

\section{Bisimilarity of diagrams up to homotopy}
\label{sec:def-bisim-up-to-homotopy)}

Let $\K$ be a category. We gather some basic results about the category $\diag(\K)$ of all small diagrams over all small categories defined as follows. An object of $\diag(\K)$ is a functor $F:{I}\to \K$ from a small category ${I}$ to $\K$. A morphism from $F:{I}_1\to \K$ to $G:{I}_2\to \K$ is a pair $(f:{I}_1\to {I}_2,\mu:F \Rightarrow G.f)$ where $f$ is a functor and $\mu$ is a natural transformation. If $(g,\nu)$ is a map from $G:{I}_2\to \K$ to $H:{I}_3\to \K$, then the composite $(g,\nu).(f,\mu)$ is defined by $(g.f,(\nu.f)\odot\mu)$ where $\odot$ means the composition of natural transformations. The identity of $F:{I}_1\to \K$ is the pair $(\id_{{I}_1},\id_F)$. It is well-known that when $\K$ is locally presentable, the category $\diag(\K)$ is locally presentable as well. I learnt the result from \cite{271951} and from a remark after the question \cite{266597}. For the convenience of the reader, the argument is recalled. First, we observe that the forgetful functor from $\diag(\K)$ to the category of small categories $\cat$ is a bifibred category by e.g. \cite[Proposition~A.1]{leftproperflow}. It corresponds to an accessible pseudo-functor in the sense of \cite[Definition~5.3.1]{MR1031717} and we use \cite[Theorem~5.3.4]{MR1031717} to deduce that $\diag(\K)$ is accessible. Finally, we observe that $\diag(\K)$ is cocomplete (e.g. \cite[Proposition~15]{dubut_PhD}) to complete the proof.

\bd \cite[Section~2.3]{dubut_bisimilarity} Two objects $F:\underline{I} \to \K$ and $G:\underline{J}\to \K$ of $\diag(\K)$ are \emph{bisimilar} if there exists a set $\mathcal{R}$ of triples $(i,\eta,j)$ called a \textit{bisimulation} where $i$ is an object of $\underline{I}$, $j$ an object of $\underline{J}$ and $\eta:F(i)\stackrel{\iso}\to G(j)$ an isomorphism of $\K$ such that the following conditions hold: 
\begin{enumerate}[leftmargin=*]
	\item For every $i\in \Obj(\underline{I})$, there exists $(i,\eta,j)\in \mathcal{R}$ and similarly for every object $j\in \Obj(\underline{J})$, there exists $(i,\eta,j)\in \mathcal{R}$.
	\item For every triple $(i,\eta,j)$ of $\mathcal{R}$ and every map $\phi:i\to i'$ of $\underline{I}$, there exists a triple $(i',\eta',j')$ of $\mathcal{R}$ and a map $\psi:j\to j'$ of $\underline{J}$ such that there is the commutative diagram of $\K$
	\[
	\begin{tikzcd}[row sep=3em, column sep=3em]
		F(i) \arrow[r,"\eta"] \arrow[d,"F(\phi)"'] & G(j) \arrow[d,dashed,"G(\psi)"]  \\
		F(i')\arrow[r,"\eta'"]& G(j')
	\end{tikzcd}
	\]
	and symmetrically, for every triple $(i,\eta,j)$ of $\mathcal{R}$ and every map $\psi:j\to j'$ of $\underline{J}$, there exists a triple $(i',\eta',j')$ of $\mathcal{R}$ and a map $\phi:i\to i'$ of $\underline{I}$ such that there is the commutative diagram of $\K$
	\[
	\begin{tikzcd}[row sep=3em, column sep=3em]
		F(i) \arrow[r,"\eta"] \arrow[d,dashed,"F(\phi)"'] & G(j) \arrow[d,"G(\psi)"]  \\
		F(i')\arrow[r,"\eta'"]& G(j')
	\end{tikzcd}
	\]	
\end{enumerate}
\ed

\bd
Let $F:\underline{I} \to \K$ and $G:\underline{J}\to \K$ be two objects of $\diag(\K)$. A map $(f,\mu):F \to G$ is \textit{open} if $f:\underline{I}\to \underline{J}$ is surjective on objects, if every map $f(i)\to j'$ lifts to a morphism $i\to j$ (in particular $f(j)=j'$) and finally if $\mu:F\Rightarrow G.f$ is a natural isomorphism.
\ed

\bp \cite[Theorem~2]{dubut_bisimilarity} \label{thm:bisimilar-open} Two objects $F:\underline{I} \to \K$ and $G:\underline{J}\to \K$ of $\diag(\K)$ are bisimilar if and only if they are related by a span of open maps. \ep

Let $\mathcal{M}$ be a model category. We denote by $\h:\mathcal{M}\to \ho(\mathcal{M})$ the canonical functor from the model category to its homotopy category. 

\bd
Let $F:\underline{I} \to \mathcal{M}$ and $G:\underline{J}\to \mathcal{M}$ be two objects of $\diag(\mathcal{M})$. A map $(f,\mu):F \to G$ is \textit{open up to homotopy} if the map $(f,\h.\mu):\h.F \to \h.G$ of $\diag(\ho(\mathcal{M}))$ is open. The diagrams $F:\underline{I} \to \mathcal{M}$ and $G:\underline{J}\to \mathcal{M}$ are \textit{bisimilar up to homotopy} if the diagram $\h.F:\underline{I} \to \ho(\mathcal{M})$ and $\h.G:\underline{J}\to \ho(\mathcal{M})$ are bisimilar.
\ed

\bp \label{prop:h-open2bisim}  
If two diagrams $F:\underline{I} \to \mathcal{M}$ and $G:\underline{J}\to \mathcal{M}$ of $\diag(\mathcal{M})$ are related by a span of maps which are open up to homotopy, then the diagrams $F:\underline{I} \to \mathcal{M}$ and $G:\underline{J}\to \mathcal{M}$ are bisimilar up to homotopy.
\ep

\bpf
Let $F:\underline{I} \to \mathcal{M}$ and $G:\underline{J}\to \mathcal{M}$ be two objects of $\diag(\mathcal{M})$ which are related by a span of maps which are open up to homotopy. This means that there exists a diagram $H$ of $\diag(\mathcal{M})$ and a span $F\leftarrow H \rightarrow G$ in $\diag(\mathcal{M})$ such that the span $\h.F\leftarrow \h.H \rightarrow \h.G$ is a span of open maps in $\diag(\ho(\mathcal{M}))$. Thus the diagrams $\h.F:\underline{I} \to \h(\mathcal{M})$ and $\h.G:\underline{J}\to \h(\mathcal{M})$ are bisimilar by Proposition~\ref{thm:bisimilar-open}. This means that the diagrams $F:\underline{I} \to \mathcal{M}$ and $G:\underline{J}\to \mathcal{M}$ are bisimilar up to homotopy. 
\epf

The converse is not true in general. Indeed, suppose that $F:\underline{I} \to \mathcal{M}$ and $G:\underline{J}\to \mathcal{M}$ are bisimilar up to homotopy. Using Proposition~\ref{thm:bisimilar-open}, we can only conclude the existence of a diagram $H:\underline{K}\to \ho(\mathcal{M})$ and of a span of open maps $\h.F\leftarrow H \rightarrow \h.G$ in $\diag(\ho(\mathcal{M}))$.

\section{Natural system associated with a general multipointed d-space}
\label{sec:construction-NT}

After \cite[Section~1]{CohomologySmallCategories}, the \textit{category of factorizations} of a small category $\C$, denoted by $\mathcal{F}(\C)$ has for objects the morphisms of $\C$ and for morphisms the \textit{extensions} of the morphisms of $\C$. This means that a morphisms from $f:X\to Y$ to $f':X'\to Y'$ is a pair of morphisms $(u:X'\to X,v:Y\to Y')$ such that there is the commutative square 
\[
\begin{tikzcd}[row sep=3em, column sep=3em]
	X \arrow[r,"f"] & Y \arrow[d,"v"] \\
    X'\arrow[u,"u"] \arrow[r,"f'"] & Y'.
\end{tikzcd}
\]
Let $\C$ be a small category. Let $\D$ be another category. A \textit{natural system} of objects of $\D$ on $\C$ is a functor $\mathcal{F}(\C) \to \D$.

\bd \label{def:trace}
Let $X$ be a directed space. The \textit{category of traces} of $X$, denoted by $\vec{T}(X)$, has for objects the points of $X$ and the set of maps $\vec{T}(X)(a,b)$ from $a\in X$ to $b\in X$ is the set of traces $\tr{\gamma}$ of directed paths $\gamma$ going from $a$ to $b$, i.e. the set of directed paths from $a$ to $b$ up to reparametrization by a map of $\mathcal{M}(1,1)$. The composition of traces, denoted by $*$, is induced by the normalized composition of directed paths, i.e. $\tr{\gamma} * \tr{\gamma'} = \tr{\gamma *_N\gamma'}$. It is strictly associative. The mapping $X\mapsto \vec{T}(X)$ induces a functor \[\vec{T}:\ptop{}\longrightarrow \cat\] from the categories of directed spaces to the one of small categories. 
\ed

Every set $\vec{T}(X)(a,b)$ is equipped with the quotient topology in $\top$ of the $\Delta$-kelleyfi\-cation of the compact-open topology by the equivalence relation induced by the identifications $\gamma\sim \gamma\phi$ for all directed paths $\gamma$ from $a$ to $b$ and all maps $\phi\in \mathcal{M}(1,1)$.

After \cite[Section~6.4]{dubut_PhD}, we associate with any directed space $X$ a natural system \[\NT(X):\mathcal{F}(\vec{T}(X)) \longrightarrow \top\] of topological spaces on $\vec{T}(X)$ as follows. The topological space $\NT(X)(\tr{\gamma})$ is by definition the topological space $\vec{T}(X)(\gamma(0),\gamma(1))$. The image of an extension of traces $\tr{\alpha * - * \beta}$ is the continuous map from $\vec{T}(X)(\gamma(0),\gamma(1))$ to $\vec{T}(X)(\alpha(0),\beta(1))$ which takes $\tr{\Gamma}$ to $\tr{\alpha}*\tr{\Gamma}*\tr{\beta}$. The mapping $X\mapsto \NT(X)$ induces a well-defined functor \[\NT:\ptop{} \longrightarrow \diag(\top)\] from the category of directed spaces to that of all small diagrams of topological spaces: see \cite[Section~6.4]{dubut_PhD} for further details.

\begin{nota}
	Using Theorem~\ref{thm:reflection}, we extend the functors $\vec{T}:\ptop{}\longrightarrow \cat$ and $\NT:\ptop{} \longrightarrow \diag(\top)$ to functors 
	\begin{align*}
		& \vec{T}:\ptop{\mathcal{M}}\longrightarrow \cat\\
		& \NT:\ptop{\mathcal{M}} \longrightarrow \diag(\top)
	\end{align*}
	by setting $\vec{T}(X)=\vec{T}(\cont(X))$ and $\NT(X)=\NT(\cont(X))$ for all multipointed $d$-spaces $X$.
\end{nota}

\section{From directed paths to traces}
\label{sec:path2trace}

This section recalls an important fact about the passage from directed paths to traces, i.e. directed paths up to reparametrization. Proposition~\ref{prop:colim} is used in the proofs of Proposition~\ref{lem:concatenation-1}, Proposition~\ref{lem:concatenation-2}, Proposition~\ref{lem:concatenation-3} and Proposition~\ref{calcul-continuous-case-1} for which the point is to verify that some obvious continuous bijections are actually homeomorphisms. It is not automatic even in the setting of $\Delta$-generated spaces since there exists a continuous bijection between Hausdorff $\Delta$-generated spaces which is a homotopy equivalence and which is not a homeomorphism: consider the $1$-dimensional sphere $\mathbf{S}^1$, the discretization $(\mathbf{S}^1)^\delta$ and the map between unreduced cones $C((\mathbf{S}^1)^\delta)\to C(\mathbf{S}^1)$ \cite{376474}.

A \textit{flow} $X$ is a small enriched semicategory. Its set of objects (preferably called \textit{states}) is denoted by $X^0$ and the space of morphisms (preferably called \textit{execution paths}) from $\alpha$ to $\beta$ is denoted by $\P_{\alpha,\beta}Y$ (e.g. \cite[Definition~10.1]{Moore1}). For any topological space $Z$, the flow $\glob(Z)$ is the flow having two states $0$ and $1$ and such that the only nonempty space of execution paths, when $Z$ is nonempty, is $\P_{0,1}\glob(Z)=Z$. 

There is a unique functor \[\dcat:\ptop{\mathcal{M}}\longrightarrow \dtop\] from the category of multipointed $d$-spaces to the category of flows taking a multipointed $d$-space $X$ to the unique flow $\dcat(X)$ such that $\dcat(X)^0=X^0$ and such that $\P_{\alpha,\beta}\dcat(X)$ is the quotient of the space of execution paths $\P_{\alpha,\beta}^{top}X$ by the equivalence relation generated by the reparametrization, the composition of $\dcat(X)$ being the composition of traces described in Definition~\ref{def:trace}.

The functor $\dcat:\ptop{\mathcal{M}}\to \dtop$ is not a left adjoint by a proof similar to the proof of \cite[Theorem~7.3]{mdtop}. Since both $\ptop{\mathcal{M}}$ and $\dtop$ are locally presentable, this implies that $\dcat:\ptop{\mathcal{M}}\to \dtop$ is not colimit-preserving. 

However, the point is that, sometimes, the functor $\dcat:\ptop{\mathcal{M}}\to \dtop$ commutes with colimits. Proposition~\ref{prop:colim} should have been put in \cite{Moore3}: it is an omission. 

\bp\label{prop:colim} Consider a pushout diagram of multipointed $d$-spaces of the form 
\[
\begin{tikzcd}[row sep=3em, column sep=3em]
	\globM(\mathbf{S}^{n-1}) \arrow[r] \arrow[d] & A \arrow[d] \\
    \globM(\mathbf{D}^{n}) \arrow[r] & \cocartesian B
\end{tikzcd}
\]
with $A$ cellular and $n\geq 0$. Then there is a pushout diagram of flows 
\[
\begin{tikzcd}[row sep=3em, column sep=3em]
	\glob(\mathbf{S}^{n-1}) \arrow[r] \arrow[d] & \dcat(A) \arrow[d] \\
\glob(\mathbf{D}^{n}) \arrow[r] & \cocartesian \dcat(B).
\end{tikzcd}
\]
\ep

\bpf[Sketch of proof]
The argument is expounded in a slightly different context in \cite[Corollary~8.12]{Moore2}. To save the reader having to read \cite{Moore2}, we outline it. The pushout diagram of multipointed $d$-spaces gives rise using \cite[Corollary~8]{Moore3} to a pushout diagram of Moore flows ($\moore^{top}$ is some \textit{right} adjoint from multipointed $d$-spaces to Moore flows)
\[
\begin{tikzcd}[row sep=3em, column sep=3em]
	\moore^{top}(\globM(\mathbf{S}^{n-1})) \arrow[r] \arrow[d] & \moore^{top}(A) \arrow[d] \\
	\moore^{top}(\globM(\mathbf{D}^n)) \arrow[r] & \cocartesian \moore^{top}(B).
\end{tikzcd}
\]
It is not necessary to recall the definitions of a Moore flow and of the functor $\moore^{top}$ from multipointed $d$-spaces to Moore flows. The only point what matters is that there exists a \textit{left} adjoint $\lmoore$ from Moore flows to flows by \cite[Theorem~10.7]{Moore1} and that, by \cite[Theorem~15]{Moore3}, there is the isomorphism of functors $\dcat \iso \lmoore\moore^{top}$. We obtain the pushout diagram of flows 
\[
\begin{tikzcd}[row sep=3em, column sep=3em]
	\dcat(\globM(\mathbf{S}^{n-1})) \arrow[r] \arrow[d] & \dcat(A) \arrow[d] \\
	\dcat(\globM(\mathbf{D}^n)) \arrow[r] & \cocartesian \dcat(B).
\end{tikzcd}
\]
It remains to observe that $\dcat(\globM(Z)) =\glob(Z)$ for any topological space $Z$ to complete the proof.
\epf

\section{Natural systems associated with a cellular multipointed d-space}
\label{sec:construction-NTd}

\begin{nota}
	In the whole section, $X$ stands for a cellular multipointed $d$-space.
\end{nota}

Let $c$ be a globular cell of $X$. The \textit{center} of $c$, denoted by $\ct{c}$, is equal to $c$ for $c\in X^0$ and is the point $\widehat{g_c}((0,\dots,0),1/2)$ when $c$ is a globular cell of dimension $\geq 1$, where $\widehat{g_c}:\globM(\mathbf{D}^{\dim(c)}) \to X$ is the attaching map of $c$.   Note that for $\alpha\in X^0$, $\ct{\alpha}$ denotes either the center of $\alpha$, i.e. $\alpha$, or the trace of the constant directed path equal to $\alpha$.

\bp
Let $c\preceq d$ be two globular cells of $X$. Then $\vec{T}(X)(\ct{c},\ct{d})$ is a singleton and the unique element is denoted by $\tr{c,d}$.
\ep

\bpf
There are two mutually exclusive cases:
\paragraph{\textbf{(1)} $c=d^-$ and $\dim(d)\geq 1$} Let $\widehat{g_d}:\globM(\mathbf{D}^{\dim(d)})\to X$ be the attaching map of $d$. The unique element of $\vec{T}(X)(\ct{c},\ct{d})$ is the trace $\tr{t\mapsto \widehat{g_d}((0,\dots,0),t/2)}$. 

\paragraph{\textbf{(2)} $c^+=d$ and $\dim(c)\geq 1$} Let $\widehat{g_c}:\globM(\mathbf{D}^{\dim(c)})\to X$ be the attaching map of $c$. The unique element of $\vec{T}(X)(\ct{c},\ct{d})$ is the trace $\tr{t\mapsto \widehat{g_c}((0,\dots,0),1/2+t/2)}$.
\epf

Let $[c_1,\dots,c_n]$ be a discrete trace of $X$, which means that $c_1\preceq \dots \preceq c_n$. Denote by $\ct{c_1,\dots,c_n}$ the trace $\ct{c_1,c_2}*\dots * \ct{c_{n-1},c_n} \in \vec{T}(X)(\ct{c_1},\ct{c_n})$. By abuse of notation, the center $\ct{c}$ of a globular cell $c$ is identified with the trace of the constant directed path equal to the center $\ct{c}$ of the globular cell $c$. Note also that \[\dt(\ct{c_1,\dots,c_n}) = [c_1,\dots,c_n].\]

\bp \label{prop:directed-path-from-globe-1}
Let $c$ be a globular cell of $X$ with $\dim(c)\geq 1$. Let $u,v\in c\cup\{c^-,c^+\}$ such that there exists a directed path from $u$ to $v$ and such that $\{u,v\} \neq \{c^-,c^+\}$. Then there is a unique trace, denoted by $\tau_{u,v}$, from $u$ to $v$ inside $c\cup\{c^-,c^+\}$. Moreover, the hypothesis $\{u,v\} \neq \{c^-,c^+\}$  cannot be removed if $\dim(c)\geq 2$.
\ep

\begin{figure}
	\def\L{7}
\begin{tikzpicture}[scale=1.3,pn/.style={circle,inner sep=0pt,minimum width=5pt,fill=black}]
	\fill[color=gray!15] (0,0) plot[smooth] coordinates {(0,0) (\L/4,1.15) (\L/2,1.5) (3*\L/4,1.15) (\L,0)} -- cycle;
	\fill[color=gray!15] (0,0) plot[smooth] coordinates {(0,0) (\L/4,-1.15) (\L/2,-1.5) (3*\L/4,-1.15) (\L,0)} -- cycle;
	\draw[dark-red][thick] plot[smooth] coordinates {(0,0) (\L/4,1.15) (\L/2,1.5) (3*\L/4,1.15) (\L,0)};
	\draw[dark-red][thick] plot[smooth] coordinates {(0,0) (\L/4,-1.15) (\L/2,-1.5) (3*\L/4,-1.15) (\L,0)};
	\draw (0,0) node[pn] {} node[black,above left]{$c^-$};
	\draw (\L/4,0.15) node[black,above]{$\overline{\tau}_x^{<}$};
	\draw (3*\L/4,0.15) node[black,above]{$\overline{\tau}_x^{>}$};
	\draw (\L/2,0.7) node[black,above]{$\overline{\tau}_x=\overline{\tau}_x^{<}*\overline{\tau}_x^{>}$};
	\draw (\L,0) node[pn] {} node[black,above right]{$c^+$};
	\draw (\L/2,0.3) node[pn] {} node[black,above]{$x$};
	\draw[-] plot[smooth] coordinates {(0,0) (\L/2,0.3) (\L,0)};
	\def\courbe {plot[smooth] coordinates {(0,0) (\L/3,-0.5) (2*\L/3,-0.5) (\L,0)}}
	\begin{scope}
		\clip (\L/3,-1) rectangle (2*\L/3,0);
		\draw[-] \courbe;
	\end{scope}
	\begin{scope}
		\clip (0,0) rectangle (\L/3,-0.5);
		\draw[dashed] \courbe;
	\end{scope}
	\begin{scope}
		\clip (2*\L/3,-0.5) rectangle (\L,0);
		\draw[dashed] \courbe;
	\end{scope}
	\draw (\L/2,-0.6) node[black,above]{$\tau_{u,v}$};
	\draw (\L/3,-0.5) node[pn] {} node[black,above]{$u$};
	\draw (2*\L/3,-0.5) node[pn] {} node[black,above]{$v$};
	\draw[->] (0,-1.8) -- (\L,-1.8) node[right]{time};
\end{tikzpicture}
	\caption{Inside the globular cell $c$}
	\label{fig:inside-c-1}
\end{figure}
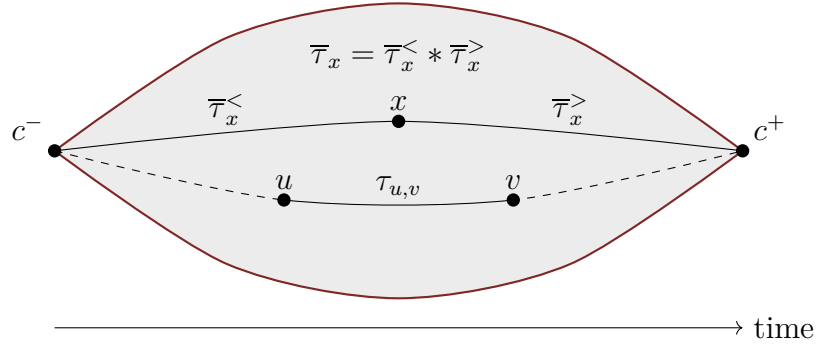

\bpf
The case $\dim(c)=1$ is obvious. We assume now $\dim(c)\geq 2$. Let \[\widehat{g}:\globM(\mathbf{D}^{\dim(c)})\longrightarrow X\] be the attaching map of the globular cell $c$. Since $\{u,v\} \neq \{c^-,c^+\}$, either $u=\widehat{g}(z,t)$ or $v=\widehat{g}(z,t)$ for some $z\in \mathbf{D}^{\dim(c)}\backslash \mathbf{S}^{\dim(c)-1}$ and some $t\in ]0,1[$. Assume the first case \[u=\widehat{g}(z,t),\] the other case being similar. Since there is a directed path from $u$ to $v$, this implies that \[v=\widehat{g}(z,t')\] for some $t\leq t'\leq 1$. The directed path from $u$ to $v$ are of the form $\widehat{g}\delta_{z}\phi$ where $\phi:[0,1]\to [t,t']$ is non-decreasing and surjective. Moreover the point $z\in \mathbf{D}^{\dim(c)}\backslash \mathbf{S}^{\dim(c)-1}$ is unique. Hence the proof of the main statement is complete. The hypothesis $\{u,v\} \neq \{c^-,c^+\}$ is necessary: if $u=c^-$ and $v=c^+$, then all execution paths $\widehat{g}\delta_{z}$ for $z$ running over $\mathbf{D}^{\dim(c)}\backslash \mathbf{S}^{\dim(c)-1}$ goes from $c^-$ to $c^+$. See Figure~\ref{fig:inside-c-1} for an illustration. 
\epf

\bp \label{prop:directed-path-from-globe-2}
Let $c$ be a globular cell of $X$. Let $x\in c$. Then there is a unique trace, denoted by $\overline{\tau}_{x}$ from $c^-$ to $c^+$ passing by $x$. One has $\overline{\tau}_{x} = \overline{\tau}_{x}^{<} * \overline{\tau}_{x}^{>}$ where $\overline{\tau}_{x}^{<}$ is the unique trace going from $c^-$ to $x$ and where $\overline{\tau}_{x}^{>}$ is the unique trace going from $x$ to $c^+$. In particular, one has $\overline{\tau}_{\ct{c}} = \overline{\tau}_{\ct{c}}^{<} =\overline{\tau}_{\ct{c}}^{>} = \tr{c}$ if $\dim(c)=0$ and $\overline{\tau}_{\ct{c}} = \tr{c^-,c,c^+}$, $\overline{\tau}_{\ct{c}}^{<} = \tr{c^-,c}$ and $\overline{\tau}_{\ct{c}}^{>} = \tr{c,c^+}$ if $\dim(c)\geq 1$.
\ep

\bpf
There is nothing to prove for $\dim(c)=0$. Assume that $\dim(c)\geq 1$. Let $\widehat{g}:\globM(\mathbf{D}^{\dim(c)})\to X$ be the attaching map of the globular cell $c$. Since $x\in c$, there exists a unique $z\in \mathbf{D}^{\dim(c)}\backslash \mathbf{S}^{\dim(c)-1}$ and a unique $t\in ]0,1[$ such that $x=\widehat{g}(z,t)$. The unique trace $\overline{\tau}_x$ is $\tr{\widehat{g}\delta_z}$. One necessarily has $\overline{\tau}_{x}^{<}=\tr{\widehat{g}\delta_z(t\mapsto (t/2))}$ and $\overline{\tau}_{x}^{>}=\tr{\widehat{g}\delta_z(t\mapsto (t/2+1/2))}$ (see Figure~\ref{fig:inside-c-1}).
\epf

\bp \label{prop:directed-path-from-globe-3}
	Let $c$ be a globular cell of $X$ with $\dim(c)\geq 1$. Let $x,y\in c$. if there exists a directed path from $x$ to $y$, then $\overline{\tau}_x = \overline{\tau}_y$. 
\ep

\bpf
Let $\widehat{g}:\globM(\mathbf{D}^{\dim(c)})\to X$ be the attaching map of the globular cell $c$. Since $x\in c$, there exists a unique $z\in \mathbf{D}^{\dim(c)}\backslash \mathbf{S}^{\dim(c)-1}$ and a unique $t\in ]0,1[$ such that $x=\widehat{g}(z,t)$. Since $y\in c$ and since there exists a directed path from $x$ to $y$, one has $x=\widehat{g}(z,t')$ for some $t\leq t'< 1$. Thus $\overline{\tau}_x = \overline{\tau}_y = \tr{\widehat{g}\delta_z}$.
\epf

\begin{lem} \label{lem:obvious}
	Let $\gamma$ be a directed path of $X$. There are five mutually exclusive cases: 
	\begin{enumerate}
		\item $\gamma(0),\gamma(1)\in X^0$;
		\item $\gamma(0)\notin X^0$, $\gamma(1)\in X^0$;
		\item $\gamma(0)\in X^0$, $\gamma(1)\notin X^0$;
		\item $\gamma(0)\notin X^0$, $\gamma(1)\notin X^0$ and $\dt(\gamma(0))\neq \dt(\gamma(1))$;
		\item $\gamma(0)\notin X^0$, $\gamma(1)\notin X^0$ and $\dt(\gamma(0))= \dt(\gamma(1))$.
	\end{enumerate}
The first four cases correspond to the situation $\dt(\gamma(0))\neq \dt(\gamma(1))$ or $\gamma(0)= \gamma(1)\in X^0$. The fifth case corresponds to the situation $\dt(\gamma(0)) = \dt(\gamma(1)) = c$ and $\dim(c)\geq 1$.
\end{lem}

\bpf Obvious.
\epf

\bp \label{lem:concatenation-1}
	Let $c$ be a globular cell of dimension greater than $1$. Let $x\in c$. Let $\alpha\in X^0$. Then there is the homeomorphism
	\begin{samepage}
		\begin{align*}
		\vec{T}(X)(c^+,\alpha) &\iso \vec{T}(X)(x,\alpha) \\
		\tau &\mapsto \overline{\tau}^{>}_{x} *\tau  
	\end{align*}
	\end{samepage}
\ep

\bpf
Consider the pushout diagram of multipointed $d$-spaces 
\[
\begin{tikzcd}[row sep=3em, column sep=3em]
\globM(\mathbf{S}^{-1}) \arrow[r,"\substack{0\mapsto \overline{x}\\1\mapsto c^+}"] \arrow[d] & \{\overline{x}\} \sqcup X \arrow[d] \\
\globM(\mathbf{D}^0) \arrow[r] & \cocartesian \overline{X}.
\end{tikzcd}
\]
The key observation is that there is the homeomorphism 
\[
\vec{P}(X)(x,\alpha) \iso \P^{top}_{\overline{x},\alpha} \overline{X}
\]
between the space of directed paths from $x$ to $\alpha$ in $X$ and the space of execution paths from $\overline{x}$ to $\alpha$ in $\overline{X}$. By Proposition~\ref{prop:colim}, we obtain the pushout diagram of flows  
\[
\begin{tikzcd}[row sep=3em, column sep=3em]
	\glob(\mathbf{S}^{-1}) \arrow[r,"\substack{0\mapsto \overline{x}\\1\mapsto c^+}"] \arrow[d] \arrow[rd,phantom,"\underline{\mathbf{A}}"{font=\tiny},pos =0.5] & \dcat(\{\overline{x}\} \sqcup X) \arrow[d] \\
\glob(\mathbf{D}^0) \arrow[r,"\widehat{g}"] & \cocartesian \dcat(\overline{X}).
\end{tikzcd}
\]
From the homeomorphism $\vec{P}(X)(x,\alpha) \iso \P^{top}_{\overline{x},\alpha} \overline{X}$, we deduce the homeomorphism 
\begin{equation}
\vec{T}(X)(x,\alpha) \iso \P_{\overline{x},\alpha} \dcat(\overline{X}) \label{eq:trace} \tag{TE}.
\end{equation}
We can now conclude. First, assume that $c^+\neq \alpha$. Then we obtain the sequence of homeomorphisms
\begin{align*}
	\vec{T}(X)(x,\alpha) & \iso \P_{\overline{x},\alpha} \dcat(\overline{X}) \\
	& \iso \P_{\overline{x},c^+} \dcat(\overline{X}) \p \P_{c^+,\alpha} \dcat(\overline{X})\\
	& \iso \P_{c^+,\alpha} \dcat(\overline{X}) \\
	& = \vec{T}(X)(c^+,\alpha),
\end{align*}
the first homeomorphism by \eqref{eq:trace}, the second homeomorphism by the commutativity of $\underline{\mathbf{A}}$ and by general results about enriched semicategories, the third homeomorphism since $\P_{\overline{x},c^+} \dcat(\overline{X})$ is a singleton, and finally the equality since $c^+\neq \alpha$. Now assume that $c^+=\alpha$. Then we obtain the sequence of homeomorphisms
\begin{align*}
	\vec{T}(X)(x,c^+) & \iso \P_{\overline{x},c^+} \dcat(\overline{X}) \\
	& \iso \{\widehat{g}(0)\} \sqcup \{\widehat{g}(0)\}\p \P_{c^+,c^+} \dcat(\overline{X})\\
	& \iso  \vec{T}(X)(c^+,c^+),
\end{align*}
the first homeomorphism by \eqref{eq:trace}, the second homeomorphism  by the commutativity of $\underline{\mathbf{A}}$ and since the execution paths from $\overline{x}$ to $c^+$ in the flow $\dcat(\overline{X})$ consists of the execution path $\widehat{g}(0)$ and all compositions $\widehat{g}(0) * \gamma$ where $\gamma$ is an execution path of $\dcat(\overline{X})$ from $c^+$ to itself, and the third homeomorphisms because the space of traces from $c^+$ to itself contains also the constant path $c^+$ which is in a distinct path-connected component because $X$ is cellular.
\epf

\begin{figure}
	\def\n{5}
	\begin{tikzpicture}[black,scale=4,pn/.style={circle,inner sep=0pt,minimum width=4pt,fill=dark-red}]
		\fill [color=gray!15] (0,0) -- (0,1) -- (1,1) -- (1,0) -- cycle;
		\draw[-] [thick] (0,1) -- (0,0);
		\draw[-] [thick] (0,0) -- (1,1);
		\foreach \n in {1,2,3,4,5,6,7,8,9,10}
		{\draw[dark-red][->][thick](0,\n/10) -- (\n/10,\n/10);}
		\draw (0,0) node[pn] {} node[black,below left] {$(0\,,0)$};
		\draw (1,1) node[black,above right] {$(1\,,1)$};
		\draw (0,1) node[black,above left] {$(0\,,1)$};
		\draw (1,0) node[black,below right] {$(1\,,0)$};
	\end{tikzpicture}
	\caption{$|U|=[0,1]\p [0,1]$, $U^0=\{0\}\p [0,1] \cup \{(x,x)\mid x\in [0,1]\}$, $\P^{top}_{(0,t),(t,t)}U= \mathcal{M}(1,1)$ for all $t\in ]0,1]$, $\P^{top}_{(0,0),(0,0)}U= \{(0,0)\}$ and $\P^{top}_{\alpha,\beta}U=\varnothing$ otherwise, there are no composable execution paths.}
	\label{fig:contracting}
\end{figure}
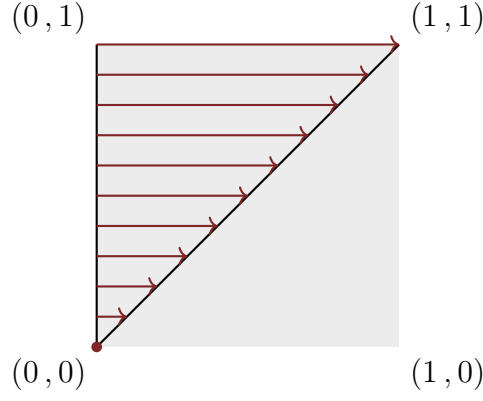

The last sentence does require for $X$ to be cellular. Figure~\ref{fig:contracting} shows an example of a multipointed $d$-space $U$ containing a sequence of non-constant execution paths converging to a constant execution path.

\bp \label{lem:concatenation-2}
	Let $c$ be a globular cell of dimension greater than $1$. Let $x\in c$. Let $\alpha\in X^0$. Then there is the homeomorphism
	\begin{samepage}
		\begin{align*}
		\vec{T}(X)(\alpha,c^-) &\iso \vec{T}(X)(\alpha,x) \\
		\tau &\mapsto  \tau * \overline{\tau}^{<}_{x}  
	\end{align*}
	\end{samepage}
\ep

\bpf We start from the pushout diagram of multipointed $d$-spaces 
\[
\begin{tikzcd}[row sep=3em, column sep=3em]
	\globM(\mathbf{S}^{-1}) \arrow[r,"\substack{0\mapsto c^-\\1\mapsto \overline{x}}"] \arrow[d] & \{\overline{x}\} \sqcup X \arrow[d] \\
\globM(\mathbf{D}^0) \arrow[r] & \cocartesian \overline{X}.
\end{tikzcd}
\]
The key observation is that there is the homeomorphism 
\[
\vec{P}(X)(\alpha,x) \iso \P^{top}_{\alpha,\overline{x}} \overline{X}
\]
between the space of directed paths from $\alpha$ to $x$ in $X$ and the space of execution paths from $\alpha$ to $\overline{x}$ in $\overline{X}$. The rest of the proof is similar to the one of Proposition~\ref{lem:concatenation-1}. It makes use of Proposition~\ref{prop:colim}.
\epf

\bp \label{lem:concatenation-3}
	Let $c\neq d$ be two globular cells of dimension greater than $1$. Let $x\in c$ and $y\in d$. Then there are the homeomorphisms 
	\begin{samepage}
		\begin{align*}
		\vec{T}(X)(c^+,d^-) &\iso \vec{T}(X)(x,y) \\
		\tau &\mapsto \overline{\tau}^{>}_{x} *\tau * \overline{\tau}^{<}_{y} 
	\end{align*}
	\end{samepage}
\ep

\bpf We start from the pushout diagram of multipointed $d$-spaces 
\[
\begin{tikzcd}[row sep=3em, column sep=3em]
	\globM(\mathbf{S}^{-1}) \sqcup \globM(\mathbf{S}^{-1})\arrow[r,"g"] \arrow[d] & \{\overline{x},\overline{y}\} \sqcup X \arrow[d] \\
\globM(\mathbf{D}^0) \sqcup \globM(\mathbf{D}^0) \arrow[r] & \cocartesian \overline{X}
\end{tikzcd}
\]
where $g$ takes the states $0,1$ of the left-hand copy of $\globM(\mathbf{S}^{-1})$ to $\overline{x}$ and $c^+$ respectively and the states $0,1$ of the right-hand copy of $\globM(\mathbf{S}^{-1})$ to $d^-$ and $\overline{y}$ respectively. The key observation is then that there is the homeomorphism 
\[
\vec{P}(X)(x,y) \iso \P^{top}_{\overline{x},\overline{y}} \overline{X}
\]
between the space of directed paths from $x$ to $y$ in $X$ and the space of execution paths from $\overline{x}$ to $\overline{y}$ in $\overline{X}$. The rest of the proof is similar to the one of Proposition~\ref{lem:concatenation-1}. It makes use of Proposition~\ref{prop:colim}.
\epf

\bp \label{calcul-continuous-case-2}
Let $\gamma$ be a directed path of $X$. Assume that $\dt(\gamma(0))\neq \dt(\gamma(1))$ or $\gamma(0)= \gamma(1)\in X^0$. Then there is the homeomorphism 
\begin{samepage}
	\begin{align*}
	\vec{T}(X)(\dt(\gamma(0))^+,\dt(\gamma(1))^-) &\iso \vec{T}(X)(\gamma(0),\gamma(1)) \\
	\tau &\mapsto \overline{\tau}^{>}_{\gamma(0)} *\tau *\overline{\tau}^{<}_{\gamma(1)} 
\end{align*}
\end{samepage}
\ep

\bpf
By Lemma~\ref{lem:obvious}, there are four mutually exclusive cases. 

\paragraph{\textbf{(1)} $\gamma(0),\gamma(1)\in X^0$} Then there is the equality \[\vec{T}(X)(\gamma(0),\gamma(1)) = \vec{T}(X)(\dt(\gamma(0))^+,\dt(\gamma(1))^-)\]
because $\gamma(0) = \dt(\gamma(0))^+ \in X^0$ and $\gamma(1) = \dt(\gamma(1))^- \in X^0$.

\paragraph{\textbf{(2)} $\gamma(0)\notin X^0$, $\gamma(1)\in X^0$} Then $\dim(\dt(\gamma(0)))\geq 1$ and one has
\begin{align*}
	\vec{T}(X)(\gamma(0),\gamma(1)) & \iso \vec{T}(X)(\dt(\gamma(0))^+,\gamma(1))\\
	 & = \vec{T}(X)(\dt(\gamma(0))^+,\dt(\gamma(1))^-),
\end{align*}
the homeomorphism by Proposition~\ref{lem:concatenation-1} and the equality because $\gamma(1)=\dt(\gamma(1))^-\in X^0$.

\paragraph{\textbf{(3)} $\gamma(0)\in X^0$, $\gamma(1)\notin X^0$} Then $\dim(\dt(\gamma(1)))\geq 1$ and one has 
\begin{samepage}
	\begin{align*}
	\vec{T}(X)(\gamma(0),\gamma(1)) & \iso \vec{T}(X)(\gamma(0),\dt(\gamma(1))^-)\\
	& = \vec{T}(X)(\dt(\gamma(0))^+,\dt(\gamma(1))^-),
\end{align*}
\end{samepage}
the homeomorphism by Proposition~\ref{lem:concatenation-2} and the equality because $\gamma(0)=\dt(\gamma(0))^+\in X^0$.

\paragraph{\textbf{(4)} $\gamma(0)\notin X^0$, $\gamma(1)\notin X^0$, $\dt(\gamma(0))\neq \dt(\gamma(1))$} Then one has $\dim(\dt(\gamma(0)))\geq 1$ and $\dim(\dt(\gamma(1)))\geq 1$. We obtain the homeomorphism
\[
\vec{T}(X)(\gamma(0),\gamma(1)) \iso \vec{T}(X)(\dt(\gamma(0))^+,\dt(\gamma(1))^-)
\]
by Proposition~\ref{lem:concatenation-3}.
\epf

\bp \label{calcul-discrete-case-2}
Let $\gamma$ be a directed path of $X$. Assume that $\dt(\gamma(0))\neq \dt(\gamma(1))$ or $\gamma(0)= \gamma(1)\in X^0$. Then there is the homeomorphism 
\begin{samepage}
	\begin{align*}
	\vec{T}(X)(\dt(\gamma(0))^+,\dt(\gamma(1))^-) &\iso \vec{T}(X)(\ct{\dt(\gamma(0))},\ct{\dt(\gamma(1))}) \\
	\tau & \mapsto \overline{\tau}^{>}_{\ct{\dt(\gamma(0))}}*\tau *\overline{\tau}^{<}_{\ct{\dt(\gamma(1)})}
\end{align*}
\end{samepage}
\ep

\bpf
We apply Proposition~\ref{calcul-continuous-case-2} to the directed path $\ct{\dt(\gamma)}$.
\epf

\bp \label{calcul-continuous-case-1}
Let $\gamma$ be a directed path of $X$. Assume that $\dt(\gamma(0))= \dt(\gamma(1))=c$ and $\dim(c)\geq 1$. Then there is the homeomorphism 
\begin{samepage}
	\begin{align*}
	\{\tau_{\gamma(0),\gamma(1)}\}\sqcup \vec{T}(X)(c^+,c^-) &\iso \vec{T}(X)(\gamma(0),\gamma(1)) \\
	\tau_{\gamma(0),\gamma(1)} &\mapsto \tau_{\gamma(0),\gamma(1)}\\
	\tau&\mapsto \overline{\tau}^{>}_{\gamma(0)} *\tau *\overline{\tau}^{<}_{\gamma(1)}
\end{align*}
\end{samepage}
\ep

\bpf
We start from the pushout diagram of multipointed $d$-spaces 
\[
\begin{tikzcd}[row sep=3em, column sep=3em]
	\globM(\mathbf{S}^{-1}) \sqcup \globM(\mathbf{S}^{-1}) \arrow[r,"g"] \arrow[d] & \{\gamma_0,\gamma_1\} \sqcup X \arrow[d] \\
\globM(\mathbf{D}^0) \sqcup \globM(\mathbf{D}^0) \arrow[r] & \cocartesian \overline{X}
\end{tikzcd}
\]
where $g$ takes the states $0,1$ of the left-hand copy of $\globM(\mathbf{S}^{-1})$ to $\gamma_0$ and $c^+$ respectively, and the states $0,1$ of the right-hand copy of $\globM(\mathbf{S}^{-1})$ to $c^-$ and $\gamma_1$ respectively. The key observation is that the space of directed paths from $\gamma(0)$ to $\gamma(1)$ in $X$ is the disjoint sum of the space of directed paths from $\gamma(0)$ to $\gamma(1)$ inside the globular cell $c$ and of the space of directed paths from $\gamma(0)$ to $\gamma(1)$ exiting the globular cell $c$ via $c^+$ and returning to the globular $c$, necessarily via $c^-$, to come back to $\gamma(1)$. We obtain the homeomorphisms 
\[
\vec{P}(X)(\gamma(0),\gamma(1)) \iso \{\gamma(0)\} \sqcup \P^{top}_{\gamma_0,\gamma_1} \overline{X}
\]
if $\gamma(0)=\gamma(1)$ and 
\[
\vec{P}(X)(\gamma(0),\gamma(1)) \iso \P^{top}_{0,1}\vI^{top} \sqcup \P^{top}_{\gamma_0,\gamma_1} \overline{X}
\]
otherwise. The rest of the proof is similar to the one of Proposition~\ref{lem:concatenation-1}. It makes use of Proposition~\ref{prop:colim}.
\epf

\bp \label{calcul-discrete-case-1}
Let $\gamma$ be a directed path of $X$. Assume that $\dt(\gamma(0))= \dt(\gamma(1))=c$ and $\dim(c)\geq 1$. Then there is the homeomorphism 
\begin{samepage}
	\begin{align*}
	\{\ct{c}\}\sqcup \vec{T}(X)(c^+,c^-) &\iso \vec{T}(X)(\ct{\dt(\gamma(0))},\ct{\dt(\gamma(1))}) \\
	\ct{c} & \mapsto \ct{c} \\
	\tau & \mapsto \ct{c,c^+} *\tau * \ct{c^-,c}
\end{align*}
\end{samepage}
\ep

\bpf
We apply Proposition~\ref{calcul-continuous-case-1} to the directed path $\ct{\dt(\gamma)}$. 
\epf

\begin{cor} \label{cor:same}
Let $\gamma$ be a directed path of $X$. There is a homeomorphism 
\[
\NT(X)(\tr{\gamma}) \iso \vec{T}(X)(\ct{\dt(\gamma(0))},\ct{\dt(\gamma(1))}).
\]
\end{cor}

\bpf 
It is a consequence of Proposition~\ref{calcul-continuous-case-2}, Proposition~\ref{calcul-discrete-case-2}, Proposition~\ref{calcul-continuous-case-1} and Proposition~\ref{calcul-discrete-case-1}.
\epf

\bd
The \textit{category of discrete traces} of $X$, denoted by $\vec{T_d}(X)$, has for objects the globular cells of $X$ and a morphism from a globular cell $c$ to a globular cell $d$ is a discrete trace $[c_1,\dots,c_n]$ with $c_1=c$ and $c_n=d$.
\ed 

We can associate with any cellular multipointed $d$-space $X$ a natural system \[\NT_d(X):\mathcal{F}(\vec{T_d}(X)) \longrightarrow \top\] of topological spaces on $\vec{T_d}(X)$ as follows. The topological spaces $\NT_d(X)([c_1,\dots,c_n])$ is by definition $\vec{T}(X)(\ct{c_1},\ct{c_n})$. The image of an extension of traces $\ct{a_1,\dots,a_m}* - * \ct{b_1,\dots,b_n}$ is the continuous map from $\vec{T}(X)(\ct{a_m},\ct{b_1})$ to $\vec{T}(X)(\ct{a_1},\ct{b_n})$ which takes $\tr{\Gamma}$ to $\ct{a_1,\dots,a_m}*\tr{\Gamma}*\ct{b_1,\dots,b_n}$. 

Theorem~\ref{thm:open-up-to-homotopy} is the globular analogue of \cite[Theorem~28 page~157]{dubut_PhD}.

\bth \label{thm:open-up-to-homotopy}
There exists a map of natural systems \[f:\NT(X)\longrightarrow \NT_d(X)\] which is open up to homotopy.
\eth

\bpf
The mapping $\dt(-)$ induces a functor from $\vec{T}(X)$ to $\vec{T_d}(X)$. The image of a point $x\in |X|$ is the unique globular cell $\dt(x)$ containing $x$ and the image of a trace $\tr{\gamma}$ is the discrete trace $\dt(\gamma)$. 

This functor induces a functor from $\mathcal{F}(\vec{T}(X))$ to $\mathcal{F}(\vec{T_d}(X))$ which takes a trace $\tr{\gamma}$ to its discrete trace $\dt(\gamma)$. The latter functor is surjective on objects since \[\dt(\ct{c_1,\dots,c_n}) = [c_1,\dots,c_n].\] We have to prove that it lifts extension of traces. 

Let $\gamma$ be a directed path of $X$ with $\dt(\gamma) = [c_1,\dots,c_n]$. Let $c_0\preceq c_1$. There are two mutually exclusive cases: (1) $\dim(c_0)=0$, $\dim(c_1)\geq 1$, $c_0=c_1^-$ (2) $\dim(c_0)\geq 1$, $\dim(c_1)= 0$, $c_0^+=c_1$.

\paragraph{\textbf{(1)} $\dim(c_0)=0$, $\dim(c_1)\geq 1$, $c_0=c_1^-$} This means that $c_0=c_1^-$ and $\gamma(0)\in c_1$. In this case, one has 
\[
\dt(\overline{\tau}^{<}_{\gamma(0)} * \tr{\gamma}) = [c_0,c_1,\dots,c_n].
\]

\paragraph{\textbf{(2)} $\dim(c_0)\geq 1$, $\dim(c_1)= 0$, $c_0^+=c_1$} This means that $\gamma(0)=c_0^+=c_1$. In this case one has 
\[
\dt(\ct{c_0,c_1} * \tr{\gamma}) = [c_0,c_1,\dots,c_n].
\]
The treatment of the extension of traces on the other side, i.e. a globular cell $c_{n+1}$ such that $c_n\preceq c_{n+1}$ is similar.

Let $\gamma$ be a directed path of $X$. Then there are the homeomorphisms between topological spaces
\[
\NT(X)(\tr{\gamma}) \iso \vec{T}(X)(\ct{\dt(\gamma(0))},\ct{\dt(\gamma(1))}) = \NT_d(X)(\dt(\gamma)),
\]
the left-hand homeomorphism by Corollary~\ref{cor:same} and the right-hand equality by definition of $\NT_d(X)$. It remains to prove that for any extension of traces $\tr{\alpha * - * \beta}$, there is a square of topological spaces which is commutative up to homotopy
\[
\begin{tikzcd}[row sep=3em, column sep=3em]
\NT(X)(\tr{\gamma}) \arrow[r,"\iso"] \arrow[d] \arrow[rd,phantom,"\underline{\mathbf{C}}"{font=\tiny},pos =0.5] & \NT_d(X)(\dt(\gamma)) \arrow[d] \\
\NT(X)(\tr{\alpha}*\tr{\gamma}*\tr{\beta}) \arrow[r,"\iso"] & \NT_d(X)(\dt(\alpha*\gamma*\beta)).
\end{tikzcd}
\]
The square $\underline{\mathbf{C}}$ can be decomposed in two squares of topological spaces $\underline{\mathbf{R}}$ and $\underline{\mathbf{S}}$ as follows:
\[
\begin{tikzcd}[row sep=3em, column sep=3em]
	\NT(X)(\tr{\gamma}) \arrow[r,"\iso"] \arrow[d] 
	\arrow[rd,phantom,"\underline{\mathbf{R}}"{font=\tiny},pos =0.5]& \NT_d(X)(\dt(\gamma)) \arrow[d] \\
\NT(X)(\tr{\alpha}*\tr{\gamma}) \arrow[r,"\iso"]\arrow[d] 
\arrow[rd,phantom,"\underline{\mathbf{S}}"{font=\tiny},pos =0.5]
& \NT_d(X)(\dt(\alpha*\gamma)) \arrow[d]\\
\NT(X)(\tr{\alpha}*\tr{\gamma}*\tr{\beta}) \arrow[r,"\iso"] & \NT_d(X)(\dt(\alpha*\gamma*\beta)).
\end{tikzcd}
\]

It suffices to prove e.g. the commutativity up to homotopy of the top one $\underline{\mathbf{R}}$. The proof of the commutativity of the bottom one $\underline{\mathbf{S}}$ is similar.

Let $\dt(\gamma) = [c_1,\dots,c_n]$. It suffices to make the proof for a directed path $\alpha$ such that $\dt(\alpha)= [c_0]$, which implies $c_0\preceq c_1$. Remember that $c_1=\dt(\gamma(0))$ and $c_n=\dt(\gamma(1))$. The sentence $c_0\preceq c_1$ means that $c_0\neq c_1$ and either $c_0=c_1^-$ or $c_0^+ = c_1$. By Lemma~\ref{lem:obvious}, there are two mutually exclusive cases: (a) $c_1 \neq c_n$ or $c_1 = c_n\in X^0$; (b) $c_1=c_n$ and $\dim(c_1)=\dim(c_n)\geq 1$. If $c_0=c_1^-$, then $\dim(c_1)\geq 1$. 

There are therefore three mutually exclusive cases : (1) $c_0=c_1^-$ and $c_1 \neq c_n$ and $\dim(c_1)\geq 1$; (2) $c_0=c_1^-$ and $c_1=c_n$ and $\dim(c_1)\geq 1$; (3) $c_0^+ = c_1$. We are going to need to subdivide the third case in two subcases (3a) $c_0^+ = c_1$ and $c_0\neq c_n$ and (3b) $c_0^+ = c_1$ and $c_0= c_n$. Remember that in all cases, $\gamma(0)\in c_1$ and $\gamma(1)\in c_n$ by Theorem~\ref{constr-dt}.

\paragraph{\textbf{(1)} $c_0=c_1^-$ and $c_1 \neq c_n$ and $\dim(c_1)\geq 1$} In this case, one has $c_0\in X^0$. One has $\NT(X)(\tr{\gamma})\iso \vec{T}(X)(c_1^+,c_n^-)$ by Proposition~\ref{calcul-continuous-case-2}. The composite map
\begin{multline*}
	\Psi:\NT(X)(\tr{\gamma}) \stackrel{\iso} \longrightarrow \NT_d(X)([c_1,\dots,c_n]) \longrightarrow \\\NT_d(X)([c_0,c_1,\dots,c_n]) \stackrel{\iso} \longrightarrow \NT(X)(\tr{\alpha}*\tr{\gamma})
\end{multline*}
takes $\Gamma\in \vec{T}(X)(c_1^+,c_n^-)$ to $\tr{c_0,c_1,c_1^+} * \Gamma$. The direct route \[\NT(X)(\tr{\gamma})\iso \vec{T}(X)(c_1^+,c_n^-) \to \NT(X)(\tr{\alpha}*\tr{\gamma}) \iso \vec{T}(X)(c_0^+,c_n^-)\] takes $\Gamma\in \vec{T}(X)(c_1^+,c_n^-)$ to $\overline{\tau}_{\gamma(0)} * \Gamma$. Since the traces $\tr{c_0,c_1,c_1^+}$ and $\overline{\tau}_{\gamma(0)}$ are the traces of two execution paths from $c_0$ to $c_1^+$  inside the globular cell $c$, they are in the same path-connected component. Thus the top square $\underline{\mathbf{R}}$ is commutative up to homotopy.

\paragraph{\textbf{(2)} $c_0=c_1^-$ and $c_1=c_n$ and $\dim(c_1)\geq 1$} In this case, one has $c_0\in X^0$. One has \[\NT(X)(\tr{\gamma})\iso \{\tau_{\gamma(0),\gamma(1)}\} \sqcup \vec{T}(X)(c_1^+,c_n^-)\] by Proposition~\ref{calcul-continuous-case-1}. Consider the composite map 
\begin{multline*}
	\Psi:\NT(X)(\tr{\gamma}) \stackrel{\iso} \longrightarrow \NT_d(X)([c_1,\dots,c_n]) \longrightarrow \\\NT_d(X)([c_0,c_1,\dots,c_n]) \stackrel{\iso} \longrightarrow \NT(X)(\tr{\alpha}*\tr{\gamma})
\end{multline*}
One has
\[
\begin{tikzcd}[row sep=0em, column sep=1em]
	\Psi:&\tau_{\gamma(0),\gamma(1)} \arrow[r,mapsto] & \ct{\tau_{\gamma(0),\gamma(1)}} =\tr{c_1} \arrow[r,mapsto] & \ct{c_0,c_1} \arrow[r,mapsto] & \ct{c_0}\\ 
\Psi:&\Gamma\in \vec{T}(X)(c_1^+,c_1^-)\arrow[r,mapsto] & \Gamma \arrow[r,mapsto] & \tr{c_0,c_1,c_1^+}*\Gamma \arrow[r,mapsto] & \tr{c_0,c_1,c_1^+}*\Gamma
\end{tikzcd}
\]

The direct route \[\NT(X)(\tr{\gamma}) \iso \{\tau_{\gamma(0),\gamma(1)}\} \sqcup \vec{T}(X)(c_1^+,c_1^-) \to \NT(X)(\tr{\alpha}*\tr{\gamma}) \iso \vec{T}(X)(c_0^+,c_1^-)\] takes $\tau_{\gamma(0),\gamma(1)}$ to $\ct{c_0}$ and $\Gamma\in \vec{T}(X)(c_1^+,c_1^-)$ to $\overline{\tau}_{\gamma(0)} * \Gamma$. Since the traces $\tr{c_0,c_1,c_1^+}$ and $\overline{\tau}_{\gamma(0)}$ are the traces of two execution paths from $c_0$ to $c_1^+$ which are in the same path-connected component, we conclude that the top square $\underline{\mathbf{R}}$ is commutative up to homotopy.

\paragraph{\textbf{(3a)} $c_0^+ = c_1$ and $c_0\neq c_n$} In this case, one has $\dim(c_0)\geq 1$ and $c_1\in X^0$. One has $\NT(X)(\tr{\gamma})\iso \vec{T}(X)(c_1^+,c_n^-)$ by Proposition~\ref{calcul-continuous-case-2}. The composite map
\begin{multline*}
	\Psi:\NT(X)(\tr{\gamma}) \stackrel{\iso} \longrightarrow \NT_d(X)([c_1,\dots,c_n]) \longrightarrow \\\NT_d(X)([c_0,c_1,\dots,c_n]) \stackrel{\iso} \longrightarrow \NT(X)(\tr{\alpha}*\tr{\gamma})
\end{multline*}
is a homeomorphism since $c_0^+=c_1^+$. The direct route \[\NT(X)(\tr{\gamma})\iso \vec{T}(X)(c_1^+,c_n^-) \stackrel{\iso}\longrightarrow \NT(X)(\tr{\alpha}*\tr{\gamma}) \iso \vec{T}(X)(c_0^+,c_n^-)\]
is a homeomorphism as well since $c_0^+=c_1^+$. We conclude that the top square $\underline{\mathbf{R}}$ is strictly commutative.

\paragraph{\textbf{(3b)} $c_0^+ = c_1$ and $c_0=c_n$} In this case, one has $\dim(c_0)\geq 1$ and $c_1\in X^0$. One still has $\NT(X)(\tr{\gamma})\iso \vec{T}(X)(c_1^+,c_n^-)$ by Proposition~\ref{calcul-continuous-case-2}. However, by Proposition~\ref{calcul-discrete-case-1}, the topological space $\NT_d(X)([c_0,c_1,\dots,c_n])$ has an additional path-connected component $\{\ct{c_0}\}$ which is not in the image of $\NT_d(X)([c_1,\dots,c_n])\to \NT_d(X)([c_0,c_1,\dots,c_n])$. Thus the composite map 
\begin{multline*}
	\Psi:\NT(X)(\tr{\gamma}) \stackrel{\iso} \longrightarrow \NT_d(X)([c_1,\dots,c_n]) \longrightarrow \\\NT_d(X)([c_0,c_1,\dots,c_n]) \stackrel{\iso} \longrightarrow \NT(X)(\tr{\alpha}*\tr{\gamma})
\end{multline*}
is not anymore a homeomorphism, unlike in \textbf{(3a)}. The direct route \[\NT(X)(\tr{\gamma})\iso \vec{T}(X)(c_1^+,c_n^-) \longrightarrow \NT(X)(\tr{\alpha}*\tr{\gamma}) \iso \vec{T}(X)(c_0^+,c_n^-)\]
is not a homeomorphism either. But the top square $\underline{\mathbf{R}}$ is still strictly commutative.
\epf

\section{Invariance under globular subdivision}
\label{sec:globular-sbd}

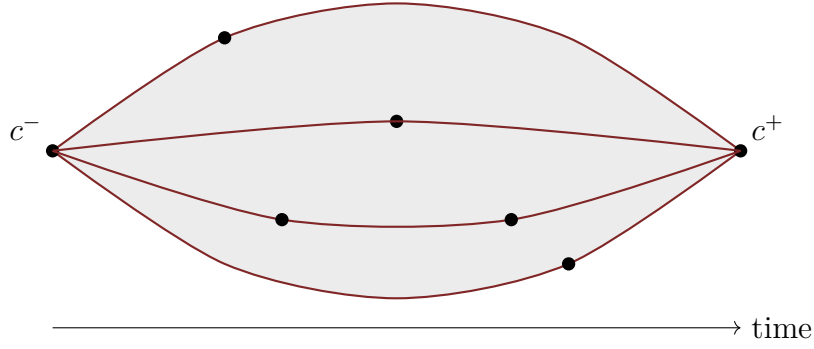
\begin{figure}
	\def\L{7}
	\begin{tikzpicture}[scale=1.3,pn/.style={circle,inner sep=0pt,minimum width=5pt,fill=black}]
		\fill[color=gray!15] (0,0) plot[smooth] coordinates {(0,0) (\L/4,1.15) (\L/2,1.5) (3*\L/4,1.15) (\L,0)} -- cycle;
		\fill[color=gray!15] (0,0) plot[smooth] coordinates {(0,0) (\L/4,-1.15) (\L/2,-1.5) (3*\L/4,-1.15) (\L,0)} -- cycle;
		\draw[dark-red][thick] plot[smooth] coordinates {(0,0) (\L/4,1.15) (\L/2,1.5) (3*\L/4,1.15) (\L,0)};
		\draw[dark-red][thick] plot[smooth] coordinates {(0,0) (\L/4,-1.15) (\L/2,-1.5) (3*\L/4,-1.15) (\L,0)};
		\draw (0,0) node[pn] {} node[black,above left]{$c^-$};
		\draw (\L,0) node[pn] {} node[black,above right]{$c^+$};
		\draw (\L/2,0.3) node[pn] {};
		\draw (\L/4,1.15) node[pn] {};
		\draw (3*\L/4,-1.15) node[pn] {};
		\draw[dark-red][thick] plot[smooth] coordinates {(0,0) (\L/2,0.3) (\L,0)};
		\def\courbe {plot[smooth] coordinates {(0,0) (\L/3,-0.7) (2*\L/3,-0.7) (\L,0)}}
		\draw[dark-red][thick] \courbe;
		\draw (\L/3,-0.7) node[pn] {};
		\draw (2*\L/3,-0.7) node[pn] {};
		\draw[->] (0,-1.8) -- (\L,-1.8) node[right]{time};
	\end{tikzpicture}
	\caption{Example of a globular subdivision of $c$}
	\label{fig:globular-sdb}
\end{figure}

\bd \cite[Definition~4.10]{diCW} \label{def:globular-sbd}
A map of multipointed $d$-spaces $f:X\to Y$ is a \textit{globular subdivision} if both $X$ and $Y$ are cellular and if $f$ induces a homeomorphism between the underlying topological spaces of $X$ and $Y$. 
\ed

The terminology used in \cite[Definition~4.10]{diCW} is \textit{T-homotopy equivalence}. The reason of the terminology of globular subdivision is that it is a globular analogue of the cubical subdivision studied in \cite{dubut_PhD}. 

\bp
There exists a map of multipointed $d$-spaces $f:X\to Y$ inducing a homeomorphism $f:|X|\to |Y|$ such that $\cont(f):\cont(X) \to \cont(Y)$ is not an isomorphism of directed spaces.
\ep

\bpf Let $X=([0,1],\{0,1\},\varnothing)$ and $Y=\vI^{top}$. Consider the map of multipointed $d$-spaces $f:X\to Y$ induced by the identity of $[0,1]$. The directed space $\cont(X)$ contains only constant paths whereas the directed space $\cont(Y)$ contains all continuous maps going from $a$ to $b$ with $a\leq b$.
\epf

\bth \label{thm:T-homotopy}
Let $f:X\to Y$ be a map of multipointed $d$-spaces between two cellular multipointed $d$-spaces $X$ and $Y$. The following two conditions are equivalent: 
\begin{enumerate}
	\item The map $f:X\to Y$ is a globular subdivision.
	\item The map of directed spaces $\cont(f):\cont(X) \to \cont(Y)$ is an isomorphism.
\end{enumerate}
\eth

\bpf 
The implication $(2)\Rightarrow(1)$ is obvious. Let us prove the implication $(1)\Rightarrow(2)$. Assume $(1)$. Since $f:|X|\to |Y|$ is one-to-one, it induces a one-to-one map from the set of directed paths $d(X)$ of $X$ to the set of directed paths $d(Y)$ of $Y$. It remains to prove that $f$ induces a surjective map from $d(X)$ to $d(Y)$. 

Consider at first the case of a directed path of $Y$ which is a regular execution path of the form $\gamma=\widehat{g'}\delta_{z'}$ with $\carrier(\gamma) = [c']$ for some globular cell $c'$ of $Y$ with corresponding attaching map $\widehat{g'}:\globM(\mathbf{D}^{\dim(c')})\to X$ and for some $z'\in \mathbf{D}^{\dim(c')}\backslash \mathbf{S}^{\dim(c')-1}$. Then $\gamma(1/2)\notin Y^0$. Thus $f^{-1}(\gamma(1/2))\notin X^0$. Consider the globular cell of $X$ of dimension greater than $1$
\[
c = \dt\bigg(f^{-1}(\gamma(1/2))\bigg)
\]
containing $f^{-1}(\gamma(1/2))$. Write $\widehat{g}:\globM(\mathbf{D}^{\dim(c)})\to X$ for the corresponding attaching map. Then 
\[
f^{-1}(\gamma(1/2)) = \widehat{g}(z,t)
\]
for some $z\in \mathbf{D}^{\dim(c)}\backslash \mathbf{S}^{\dim(c)-1}$ and some $t\in ]0,1[$. The execution path $f\widehat{g}\delta_z$ of $Y$ contains $\gamma(1/2)$. One has $\gamma(0) = \widehat{g'}\delta_{z'}(0) = (f\widehat{g}\delta_z)(t_0) \in Y^0$ and $\gamma(1) = \widehat{g'}\delta_{z'}(1) = (f\widehat{g}\delta_z)(t_1) \in Y^0$ for some real numbers $t_0,t_1$ such that $0\leq t_0<t<t_1\leq 1$. Consider the directed path $\Gamma$ of $Y$ defined by 
\[
\Gamma(u) = f(\widehat{g}(z,t_0+(t_1-t_0)u)).
\]
It goes from $\Gamma(0)=f(\widehat{g}(z,t_0)) = \gamma(0)$ to $\Gamma(1)=f(\widehat{g}(z,t_1)) = \gamma(1)$. Moreover, one has
\[
\Gamma\bigg(\frac{t-t_0}{t_1-t_0}\bigg) = f(\widehat{g}(z,t)) = \gamma(1/2).
\]
The directed path $\Gamma$ being regular since $f$ is a bijective, there is the inclusion \[\Gamma(]0,1[) \subset \dt(\gamma(1/2))=c'.\] This means that $\Gamma$ is a directed path from $\gamma(0)$ to $\gamma(1)$ inside the globular cell $c'$ of $Y$ passing by $\gamma(1/2)$. This implies that $\Gamma$ is an execution path of $Y$ and that there exists $\phi\in \mathcal{M}(1,1)$ such that 
\[
\Gamma = \widehat{g'}\delta_{z'} \psi' = \gamma \phi.
\]
Since both $\Gamma$ and $\gamma$ are regular paths, $\phi$ is actually a homeomorphism from $[0,1]$ to itself by Proposition~\ref{prop:regular_reparam}. We obtain 
\[
\gamma = f\bigg( \widehat{g} \delta_z \big(u\mapsto t_0+(t_1-t_0)\phi^{-1}(u)\big)\bigg)
\]
Since the map $u\mapsto t_0+(t_1-t_0)\phi^{-1}(u)\in \mathcal{I}(1)$, we have proved that $\gamma$ is the image by $f$ of a directed path of $X$. 

Using Theorem~\ref{thm:decomposition_naturelle}, we deduce that any execution path of $Y$ is the image by $f$ of a directed path of $X$. Finally, by Theorem~\ref{continuous}, we obtain that any directed path of $Y$ is the image by $f$ of a directed path of $X$. 
\epf

From Figure~\ref{fig:globular-sdb}, it seems that for any globular subdivision $f:X\to Y$, any $\alpha\in Y^0\backslash X^0$ is on a directed path between two states of $f(X^0)$ and that this directed path is locally unique up to reparametrization. Indeed, for such a state $\alpha$, the point $f^{-1}(\alpha)$ belongs to a unique globular cell $\dt(f^{-1}(\alpha))$ of $X$. There is unique trace $\overline{\tau}_{f^{-1}(\alpha)}$ going from $\dt(f^{-1}(\alpha))^-$ to $\dt(f^{-1}(\alpha))^+$. And $\alpha$ belongs to $f(\dt(f^{-1}(\alpha)))$.

\bth \label{thm:final}
	Let $f:X\to Y$ be a globular subdivision. Then the natural systems $\NT_d(X)$ and $\NT_d(Y)$ are bisimilar up to homotopy.
\eth

\bpf
It is a consequence of Proposition~\ref{prop:h-open2bisim}, Theorem~\ref{thm:open-up-to-homotopy} and Theorem~\ref{thm:T-homotopy}.
\epf

\section{Natural systems and q-model structure}
\label{sec:counter-example}

A model category is a bicomplete category $\mathcal{M}$ equipped with a class of cofibrations $\C$, a class of fibrations $\F$ and a class of weak equivalences $\W$ such that: 1) $\W$ is closed under retract and satisfies \ttt, 2) the pairs $(\C,\W\cap \F)$ and $(\C\cap \W,\F)$ are functorial weak factorization systems. We refer to \cite[Chapter~1]{MR99h:55031} and to \cite[Chapter~7]{ref_model2} for the basic notions about general model categories.

The purpose of this section is to show the incompatibility of the notion of bisimilar natural systems and the model structures studied so far on multipointed $d$-spaces. 

The \textit{q-model structure} of multipointed $d$-spaces is the unique combinatorial model structure such that 
\[\{\globM(\mathbf{S}^{n-1})\subset \globM(\mathbf{D}^{n}) \mid n\geq 0\} \cup \{C:\varnothing \to \{0\},R:\{0,1\} \to \{0\}\}\]
is the set of generating cofibrations, the maps between globes being induced by the closed inclusions $\mathbf{S}^{n-1}\subset \mathbf{D}^{n}$, and such that 
\[
\{\globM(\mathbf{D}^{n})\subset \globM(\mathbf{D}^{n+1}) \mid n\geq 0\}
\]
is the set of generating trivial cofibrations, the maps between globes being induced by the closed inclusions $(x_1,\dots,x_n)\mapsto (x_1,\dots,x_n,0)$. The weak equivalences are the maps of multipointed $d$-spaces $f:X\to Y$  inducing a bijection $f^0:X^0\iso Y^0$ and a weak homotopy equivalence $\P^{{top}} f:\P^{{top}} X \to \P^{{top}} Y$ and the fibrations are the maps of multipointed $d$-spaces $f:X\to Y$  inducing a q-fibration $\P^{{top}} f:\P^{{top}} X \to \P^{{top}} Y$ of topological spaces. 

A construction of this model structure is given in \cite[Theorem~6.16]{QHMmodel} for a more specific version of multipointed $d$-spaces. The argument still works here because it relies on the use of the Quillen path object argument in \cite[Theorem~6.14]{QHMmodel} applied to the right adjoint from multipointed $d$-spaces to topological graphs which forgets the composition and the reparametrization of execution paths.

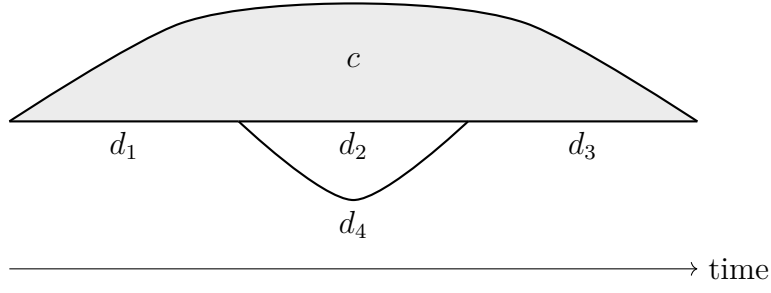
\begin{figure}
	\def\L{7}
	\begin{tikzpicture}[scale=1.3,pn/.style={circle,inner sep=0pt,minimum width=5pt,fill=black}]
		\fill[color=gray!15] (0,0) plot[smooth] coordinates {(0,0) (\L/4,1) (\L/2,1.2) (3*\L/4,1) (\L,0)} -- cycle;
		\draw[thick][-] (0,0) -- (\L/3,0);
		\draw (\L/6,0) node[below]{$d_1$};
		\draw[thick][-] (\L/3,0) -- (2*\L/3,0);
		\draw[thick][-] plot[smooth] coordinates {(\L/3,0) (1.5*\L/3,-0.8) (2*\L/3,0)};
		\draw (\L/3+\L/6,0) node[below]{$d_2$};
		\draw[thick][-] (2*\L/3,0) -- (\L,0);
		\draw (2*\L/3+\L/6,0) node[below]{$d_3$};
		\draw (1.5*\L/3,-0.8) node[below]{$d_4$};
		\draw[thick][-] plot[smooth] coordinates {(0,0) (\L/4,1) (\L/2,1.2) (3*\L/4,1) (\L,0)};
		\draw (\L/2,0.8) node[below]{$c$};
		\draw[->] (0,-1.5) -- (\L,-1.5) node[right]{time};
	\end{tikzpicture}
	\caption{The multipointed $d$-spaces $A$}
	\label{fig:A}
\end{figure}

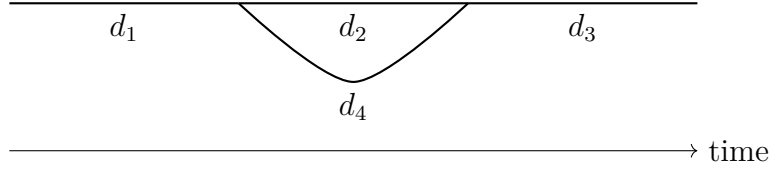
\begin{figure}
	\def\L{7}
	\begin{tikzpicture}[scale=1.3,pn/.style={circle,inner sep=0pt,minimum width=5pt,fill=black}]
		\draw[thick][-] (0,0) -- (\L/3,0);
		\draw (\L/6,0) node[below]{$d_1$};
		\draw[thick][-] (\L/3,0) -- (2*\L/3,0);
		\draw[thick][-] plot[smooth] coordinates {(\L/3,0) (1.5*\L/3,-0.8) (2*\L/3,0)};
		\draw (\L/3+\L/6,0) node[below]{$d_2$};
		\draw[thick][-] (2*\L/3,0) -- (\L,0);
		\draw (2*\L/3+\L/6,0) node[below]{$d_3$};
		\draw (1.5*\L/3,-0.8) node[below]{$d_4$};
		\draw[->] (0,-1.5) -- (\L,-1.5) node[right]{time};
	\end{tikzpicture}
	\caption{The multipointed $d$-space $B$}
	\label{fig:B}
\end{figure}

\bp \label{prop:counterexample}
There exists a trivial q-fibration $f:A\to B$ between two cellular multipointed $d$-spaces such that the map $\vec{NT}(f):\vec{NT}(A) \to \vec{NT}(B)$ is not open up to homotopy. 
\ep

\bpf
Consider the multipointed $d$-space $B$ depicted in Figure~\ref{fig:B}. Consider the multipointed $d$-space $A$, depicted in Figure~\ref{fig:A}, and obtained by the pushout diagram which adds a globular cell $c$ of dimension $2$:
\[
\begin{tikzcd}[row sep=4em, column sep=6em]
	\globM(\{0\}) \arrow[d,"{\{0\}\subset [0,1]}"'] \arrow[r,"0\mapsto d_1*_N d_2*_N d_3"] & B\arrow[d] \\
    \globM([0,1]) \arrow[r,"\widehat{g}"] & \cocartesian A
\end{tikzcd}
\]
There is a map of multipointed $d$-space $f:A\to B$ which is a trivial q-fibration and which intuitively crushes the $2$-dimensional globe $c$ depicted in Figure~\ref{fig:A}. The point is that $\NT(A)(\tr{c})$ is contractible whereas $\NT(B)(\tr{d_1,d_1^+,d_2,d_2^+,d_3})$ has two distinct path-connected components. Therefore the map $\NT(A)(\tr{c}) \to \NT(B)(\tr{d_1,d_1^+,d_2,d_2^+,d_3})$ is not a weak homotopy equivalence.
\epf

\begin{figure}
	\def\L{7}
	\begin{tikzpicture}[scale=1.3,pn/.style={circle,inner sep=0pt,minimum width=5pt,fill=black}]
		\fill[color=gray!15] (0,0) -- (\L,0) -- (\L,2) -- (0,2) -- cycle;
		\draw[thick][-] (0,0) -- (\L/3,0);
		\draw (\L/6,0) node[below]{$d_1$};
		\draw[thick][-] (\L/3,0) -- (2*\L/3,0);
		\draw[thick][-] plot[smooth] coordinates {(\L/3,0) (1.5*\L/3,-0.8) (2*\L/3,0)};
		\draw (\L/3+\L/6,0) node[below]{$d_2$};
		\draw[thick][-] (2*\L/3,0) -- (\L,0);
		\draw (2*\L/3+\L/6,0) node[below]{$d_3$};
		\draw (1.5*\L/3,-0.8) node[below]{$d_4$};
		\draw[->] (0,-1.5) -- (\L,-1.5) node[right]{time};
		\draw (1.5*\L/6,1) node[pn] {} node[black,above]{$u$};
		\draw (1.5*\L/6,0) node[pn] {} node[black,below]{$u'$};
		\draw (4.5*\L/6,1) node[pn] {} node[black,above]{$v$};
		\draw (4.5*\L/6,0) node[pn] {} node[black,below]{$v'$};
		\draw[-,dark-red] (0,1) -- (\L,1);
		\draw[dashed] (1.5*\L/6,1) -- (1.5*\L/6,0);
		\draw[dashed] (4.5*\L/6,1) -- (4.5*\L/6,0);
	\end{tikzpicture}
	\caption{The directed space $X$: the directed paths are horizontal and non-decreasing with respect to the time line}
	\label{fig:X}
\end{figure}
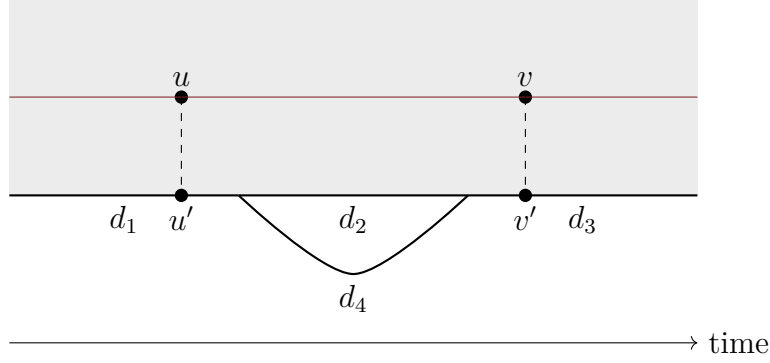

\begin{figure}
	\def\L{7}
	\begin{tikzpicture}[scale=1.3,pn/.style={circle,inner sep=0pt,minimum width=5pt,fill=black}]
		\fill[color=gray!15] (0,0) -- (\L,0) -- (\L,2) -- (0,2) -- cycle;
		\draw[thick][-] (0,0) -- (\L/3,0);
		\draw (\L/6,0) node[below]{$d_1$};
		\draw[thick][-] (\L/3,0) -- (2*\L/3,0);
		\draw[thick][-] plot[smooth] coordinates {(\L/3,0) (1.5*\L/3,-0.8) (2*\L/3,0)};
		\draw (\L/3+\L/6,0) node[below]{$d_2$};
		\draw[thick][-] (2*\L/3,0) -- (\L,0);
		\draw (2*\L/3+\L/6,0) node[below]{$d_3$};
		\draw (1.5*\L/3,-0.8) node[below]{$d_4$};
		\draw[->] (0,-1.5) -- (\L,-1.5) node[right]{time};
		\draw (1.5*\L/6,1) node[pn] {} node[black,above]{$u$};
		\draw (1.5*\L/6,0) node[pn] {} node[black,below]{$u'$};
		\draw (4.5*\L/6,1) node[pn] {} node[black,above]{$v$};
		\draw (4.5*\L/6,0) node[pn] {} node[black,below]{$v'$};
		\draw[-,dark-red] (0,1) -- (1.5*\L/6,1) -- (\L/3,0) -- (2*\L/3,0) -- (4.5*\L/6,1) -- (\L,1);
		\draw[dashed] (1.5*\L/6,1) -- (1.5*\L/6,0);
		\draw[dashed] (4.5*\L/6,1) -- (4.5*\L/6,0);
	\end{tikzpicture}
	\caption{The directed space $X'$: the directed paths are non-decreasing with respect to the time line but not necessarily horizontal}
	\label{fig:X'}
\end{figure}
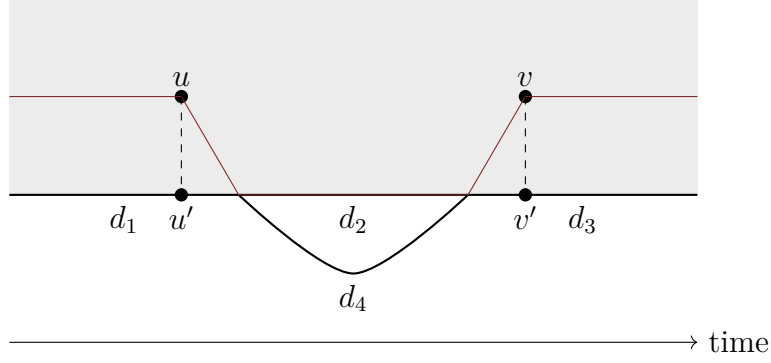

We want to conclude this section by reproducing the same geometric phenomenon with a directed space without globes. Consider the directed space $X$ depicted in Figure~\ref{fig:X} which contains $\cont(B)$ and the gray rectangle above $\cont(B)$. The directed paths of $X$ in the gray rectangle are by definition all horizontal paths which are non-decreasing with respect to the direction of time. The retract $r:X\to \cont(B)$ of the gray rectangle on the line $d_1d_2d_3$ does not induce an open map $\NT(r):\NT(X)\to \NT(B)$ since $\vec{T}(u,v)$ is contractible whereas $\vec{T}(u',v')$ contains two path-connected components. The situation is different with the directed space $X'$ depicted in Figure~\ref{fig:X'} having the same underlying topological space as $X$ and such that the directed paths in the gray rectangle are all paths (not only the horizontal ones) which are non-decreasing with respect to the direction of time. In that case, the map $\NT(r'):\NT(X')\to \NT(B)$ is an open map.

\section{Future works}
\label{sec:future}

The reason of the incompatibility, as expounded in Proposition~\ref{prop:counterexample}, is that the formalism of multipointed $d$-space in the sense of Definition~\ref{def:multipointed-d-space} and its associated q-model structure cannot see what happens between the states of a multipointed $d$-space. Multipointed $d$-spaces are a \textit{multipointed} formalism of DAT, like flows \cite{model3,Moore1,QHMmodel}, and also like the combinatorial models of precubical sets \cite{DAT_book} and of transverse sets, the latter adding transverse degeneracies to precubical sets \cite{DirectedDegeneracy,ThickCubes}. All these formalisms have indeed a distinguished set of states or vertices. On the contrary, Grandis' directed spaces \cite{mg} (see Definition~\ref{def:directed_space}) and Krishnan's streams \cite{MR2545830} are \textit{continuous} models of DAT. These two continuous models are in fact essentially the same by \cite{Ziemiaski2012}. Bisimilar natural systems preserve the causal structure even between the states of a multipointed $d$-space because this notion of bisimilarity is initially designed for the continuous model of Grandis' directed space \cite{dubut_PhD}.  

The right mathematical object to link these two approaches of DAT (namely the \textit{multipointed} one on one hand and the \textit{continuous} one on the other hand) is the categorical localization of the full subcategory of \textit{cellular} multipointed $d$-spaces in the sense of Section~\ref{sec:directed-paths-cellular-case} by the \textit{globular subdivisions} in the sense of Definition~\ref{def:globular-sbd}. This categorical localization is locally small by a variant of \cite[Theorem~4.12]{diCW}. 

All reasonable realizations of a precubical set as a cellular flow (i.e. a cellular object of the q-model structure of flows) are homotopy equivalent by \cite[Theorem~3.8]{NaturalRealization}. By \cite[Theorem~15]{Moore3}, all reasonable realizations of a precubical set as a cellular multipointed $d$-space are thus homotopy equivalent as well. However, Proposition~\ref{prop:counterexample} tells us that they are not bisimilar. Only some of them have the correct bisimilarity type. This problem will be studied in a subsequent paper.

\appendix

\section{The maximal directed subspace of a multipointed d-space}
\label{sec:opcont}

The purpose of this section is to prove that the functor $\discont:\ptop{}\to \ptop{\mathcal{M}}$ has a right adjoint which is the functor which takes a multipointed $d$-space $Y$ to its maximal directed subspace $\opcont(Y)$. At first, Proposition~\ref{prop:const-opcont} expounds the construction of the maximal directed subspace of a multipointed $d$-space.

\bp \label{prop:const-opcont}
Let $Y$ be a multipointed $d$-space. Consider the topological space \[|\opcont(Y)|=\big\{\alpha\in Y^0\mid \alpha\in \P^{top}_{\alpha,\alpha}Y\big\}\] equipped with the $\Delta$-kelleyfication of the relative topology induced by the inclusion \[|\opcont(Y)|\subset |Y|\] and the set of continuous paths \[d(\opcont(Y)) = \big\{\gamma:[0,1]\to |\opcont(Y)|\mid \forall \phi\in \mathcal{I}(1),\gamma\phi\in \P^{top}Y\big\}.\]
The pair $\opcont(Y)=(|\opcont(Y)|,d(\opcont(Y)))$ is a directed space. The mapping $Y\mapsto\opcont(Y)$ yields a well-defined functor from the category of multipointed $d$-spaces to the category of directed spaces. Finally, for any directed space $X$, a map of multipointed $d$-spaces $\discont(X)\to Y$ factors uniquely as a composite \[\discont(X)\longrightarrow \discont(\opcont(Y))\subset Y.\]
\ep

Remember that for a multipointed $d$-space $Y$, the set $Y^0$ is not necessarily discrete: see Figure~\ref{fig:contracting} for an example.

\bpf
By definition of the topology of $|\opcont(Y)|$, a set map $f:Z\to |\opcont(Y)|$ is continuous if and only if the composite set map $f:Z\to |\opcont(Y)| \subset |Y|$ is continuous. Since the identity of $[0,1]$ belongs to $\mathcal{I}(1)$, all directed paths of $\opcont(Y)$ are execution paths of $Y$.

\paragraph{\textbf{(1)}} Let $\gamma\in d(\opcont(Y))$ and $\psi\in \mathcal{I}(1)$. Then for all $\phi \in \mathcal{I}(1)$, $(\gamma\psi)\phi =\gamma(\psi\phi)\in \P^{top}Y$ since $\psi\phi\in \mathcal{I}(1)$. Thus $\gamma\psi\in d(\opcont(Y))$. 

\paragraph{\textbf{(2)}} Let $\gamma_1,\gamma_2\in d(\opcont(Y))$ such that $\gamma_1 *_N \gamma_2$ exists. Let $\phi\in \mathcal{I}(1)$. We want to prove that $(\gamma_1 *_N \gamma_2)\phi \in \P^{top}Y$. There are three possible cases (the first two are not mutually exclusive, both containing the case $\phi([0,1])=\{1/2\}$):
\begin{enumerate}[(2a)]
	\item $\phi([0,1]) \subset [0,1/2]$;
	\item $\phi([0,1]) \subset [1/2,1]$;
	\item $\phi([0,u])\subset [0,1/2]$, $\phi([u,1])\subset [1/2,1]$ with $\phi(0)<\phi(u)=1/2<\phi(1)$.
\end{enumerate}
\paragraph{\textbf{(2a)}} Then $(\gamma_1 *_N \gamma_2)\phi(t) = \gamma_1 (2\phi(t))$ for all $t\in[0,1]$. Since $t\mapsto 2\phi(t) \in \mathcal{I}(1)$, we deduce that $(\gamma_1 *_N \gamma_2)\phi\in \P^{top}Y$, $\gamma_1$ being an element of $d(\opcont(Y))$ by hypothesis.

\paragraph{\textbf{(2b)}} Then $(\gamma_1 *_N \gamma_2)\phi(t) = \gamma_2(2\phi(t)-1)$ for all $t\in[0,1]$. Since $t\mapsto 2\phi(t)-1 \in \mathcal{I}(1)$, we deduce that $(\gamma_1 *_N \gamma_2)\phi\in \P^{top}Y$, $\gamma_2$ being an element of $d(\opcont(Y))$ by hypothesis.

\paragraph{\textbf{(2c)}} This implies $0<u<1$, the map $\phi$ being non-decreasing. Write 
\[
(\gamma_1 *_N \gamma_2)\phi= \big((\gamma_1\mu_{1/2}) * (\gamma_2\mu_{1/2})\big)\phi= (\gamma_1 * \gamma_2) \underbrace{(\mu_{1/2}\ot \mu_{1/2})\phi}_{\psi}
\]
the left-hand equality by definition of the normalized composition and the right-hand equality by Proposition~\ref{lem-1}. One has $\psi\in \mathcal{I}(2)$, $\mu_{1/2}\ot \mu_{1/2}$ being a non-decreasing homeomorphism from $[0,1]$ to $[0,2]$. We deduce, by definition of the Moore composition, that (recall that $\psi(u) = (\mu_{1/2}\ot \mu_{1/2})(1/2)=1$)
\begin{equation}
	(\gamma_1 *_N \gamma_2)\phi(t) = \begin{cases}
		\gamma_1(\psi(t)) & \hbox{ if } 0\leq t\leq u\\
		\gamma_2(\psi(t)-1) & \hbox{ if } u\leq t\leq 1.
	\end{cases}
	\label{eq:moore} \tag{MO}
\end{equation}
Consider the continuous maps $\psi_1$ and $\psi_2$ defined as follows:
\begin{samepage}
	\begin{align*}
	\psi_1:[0,u] & \longrightarrow [0,1] & \psi_2:[u,1] & \longrightarrow [0,1]\\
	t & \mapsto \psi(t) & t & \mapsto \psi(t)-1
\end{align*}
\end{samepage}
One obtains the equalities (remember that $0<u<1$)
\[
(\gamma_1 *_N \gamma_2)\phi = (\gamma_1\psi_1) * (\gamma_2\psi_2) = (\gamma_1\underbrace{\psi_1\mu^{-1}_u}_{\in \mathcal{I}(1)}\mu_u) * (\gamma_2\underbrace{\psi_2\mu^{-1}_{1-u}}_{\in \mathcal{I}(1)}\mu_{1-u})
\]
the left-hand equality by \eqref{eq:moore} and the right-hand equality by trivial algebraic substitutions and because $0<u<1$ by hypothesis. Since $\gamma_1,\gamma_2\in d(\opcont(Y))$ by hypothesis, we obtain $\gamma_1\psi_1\mu^{-1}_u,\gamma_2\psi_2\mu^{-1}_{1-u} \in \P^{top}Y$. We deduce that $\gamma_1\psi_1\in \P^uY$ and $\gamma_2\psi_2\in \P^{1-u}Y$ (cf. Definition~\ref{def:length-path}). Using Proposition~\ref{prop:lengh-path}, we deduce that $(\gamma_1 *_N \gamma_2)\phi\in \P^{u+1-u}Y = \P^{top}Y$. We have proved that $\gamma_1 *_N \gamma_2$ belongs to $d(\opcont(Y))$.

This means that the pair $\opcont(Y)=(|\opcont(Y)|,d(\opcont(Y)))$ yields a well-defined directed space. Since a map of multipointed $d$-spaces $f:Y_1\to Y_2$ takes an element of $Y_1^0$ to an element of $Y_2^0$ and a constant execution path $\alpha\in \P^{top}_{\alpha,\alpha}Y_1$ to the constant execution path $f(\alpha)\in \P^{top}_{f(\alpha),f(\alpha)}Y_2$, we have obtained a functor $Y\mapsto |\opcont(Y)|$ from multipointed $d$-spaces to topological spaces. Let $\gamma\in d(\opcont(Y_1))$. Then $f\gamma\in \P^{top} Y_2$ since $f:Y_1\to Y_2$ is, by hypothesis, a map of multipointed $d$-spaces. Let $\phi\in \mathcal{I}(1)$. Then $(f\gamma)\phi = f(\gamma\phi)$ belongs to $\P^{top} Y_2$ since $\gamma\phi\in\P^{top} Y_1$, $\gamma$ being an element of $d(\opcont(Y_1))$. We have proved that the mapping $Y\mapsto \opcont(Y)$ yields a functor from the category of multipointed $d$-spaces to that of directed spaces. 

The inclusion of multipointed $d$-spaces $\discont(\opcont(Y))\subset Y$ comes from the fact that $|\discont(\opcont(Y))| = |\opcont(Y)|\subset |Y|$, $\discont(\opcont(Y))^0=|\opcont(Y)|\subset Y^0$ and from the fact that every directed path of $\opcont(Y)$ is an execution path of $Y$ by construction. 

Consider a map of multipointed $d$-spaces $f:\discont(X)\to Y$ for some directed space $X$. Then $f$ induces a set map $\discont(X)^0=|X| \to Y^0$ and a continuous composite map $|X|=|\discont(X)| \to Y^0\subset |Y|$ by construction of $\discont$ (see Proposition~\ref{prop:Omega}). For every $\alpha\in |X|$, $\alpha$ is a constant execution path of $\discont(X)$, which implies that $f(\alpha)$ is a constant execution path of $\P_{f(\alpha),f(\alpha)}Y$. The space $|\opcont(Y)|$ being equipped with the $\Delta$-kelleyfication of the relative topology, this implies that $f$ induces a continuous map from $|X|$ to $|\opcont(Y)|$. By definition of $d(\opcont(Y))$, a directed path of $X$ is mapped to a directed path of $d(\opcont(Y))$. Thus the map of multipointed $d$-spaces $f:\discont(X)\to Y$ factors as a composite $\discont(X)\longrightarrow \discont(\opcont(Y))\subset Y$. This factorization is unique because of the inclusion $|\opcont(Y)|\subset |Y|$.
\epf

The definition of the q-model structure of multipointed $d$-spaces is recalled in Section~\ref{sec:counter-example}.

\bp
	For any q-cofibrant multipointed $d$-space $X$, in particular for any cellular multipointed $d$-space $X$, one has $\opcont(X)=\varnothing$.
\ep

\bpf
A cellular multipointed $d$-space $X$ does not contain any constant execution path. This implies $\opcont(X)=\varnothing$. The proof is complete because every q-cofibrant multipointed $d$-space is a retract of a cellular one.
\epf

\bth \label{thm:opcont}
There is an adjunction $\discont\dashv \opcont$. 
\eth

\bpf
There are the natural bijections of sets 
\[
\ptop{\mathcal{M}}(\discont(X),Y) \iso \ptop{\mathcal{M}}(\discont(X),\discont(\opcont(Y))) \iso \ptop{}(X,\opcont(Y)),
\]
the left-hand bijection by Proposition~\ref{prop:const-opcont}, and the right-hand bijection by Proposition~\ref{prop:Omega}.
\epf

We conclude with the following theorem which is a complement to Theorem~\ref{thm:reflection}.

\bth \label{thm:coreflection} 
For all directed spaces $X$, there is the equality $X=\opcont(\discont(X))$. The functor $\discont:\ptop{}\to \ptop{\mathcal{M}}$ induces an isomorphism between the category of directed spaces and a full coreflective subcategory of the category of multipointed $d$-spaces. 
\eth

\bpf
The unit $X\to\opcont(\discont(X))$ of the adjunction $\discont \dashv \opcont$ induces the identity on the underlying space $|X|$. This implies that the identity of $|X|$ yields a one-to-one set map from the directed paths of $X$ to the directed paths of $\opcont(\discont(X))$. Every directed path of $\opcont(\discont(X))$ is an execution path of the multipointed $d$-space $\discont(X)$ by definition of $\opcont$, and therefore a directed path of $X$ by definition of $\discont$. We obtain the equality $X=\opcont(\discont(X))$. Finally, for every multipointed $d$-space of the form $\discont(X)$, one has $(\discont\opcont)(\discont(X)) = \discont(X)$ by the previous equality. Hence the isomorphism of categories.
\epf


\begin{thebibliography}{10}
	
	\bibitem{TheBook}
	J.~Ad{\'a}mek and J.~Rosick{\'y}.
	\newblock {\em Locally presentable and accessible categories}.
	\newblock Cambridge University Press, Cambridge, 1994.
	\newblock \href {https://doi.org/10.1017/cbo9780511600579.004}
	{\path{https://doi.org/10.1017/cbo9780511600579.004}}.
	
	\bibitem{CohomologySmallCategories}
	H.-J. Baues and G.~Wirsching.
	\newblock Cohomology of small categories.
	\newblock {\em J. Pure Appl. Algebra}, 38:187--211, 1985.
	\newblock \href {https://doi.org/10.1016/0022-4049(85)90008-8}
	{\path{https://doi.org/10.1016/0022-4049(85)90008-8}}.
	
	\bibitem{271951}
	T.~Campion.
	\newblock Presentable small diagrams over a locally presentable category.
	\newblock MathOverflow, 2017.
	\newblock URL:https://mathoverflow.net/q/271951 (version: 2017-06-11).
	\newblock \href {http://arxiv.org/abs/https://mathoverflow.net/q/271951}
	{\path{https://mathoverflow.net/q/271951}}.
	
	\bibitem{dubut_PhD}
	J.~Dubut.
	\newblock {\em Directed homotopy and homology theories for geometric models of
		true concurrency}.
	\newblock PhD thesis, Universit\'e de Saclay, Ecole Normale Sup\'erieure de
	Cachan, 2017.
	
	\bibitem{dubut_bisimilarity}
	J.~Dubut.
	\newblock Bisimilarity of diagrams.
	\newblock In {\em Relational and algebraic methods in computer science. 18th
		international conference, RAMiCS 2020, Palaiseau, France, October 26--29,
		2020. Proceedings}, pages 65--81. Cham: Springer, 2020.
	\newblock \href {https://doi.org/10.1007/978-3-030-43520-2_5}
	{\path{https://doi.org/10.1007/978-3-030-43520-2_5}}.
	
	\bibitem{NaturalHomology}
	J.~Dubut, E.~Goubault, and J.~Goubault-Larrecq.
	\newblock Natural homology.
	\newblock In {\em Automata, languages, and programming. 42nd international
		colloquium, ICALP 2015, Kyoto, Japan, July 6--10, 2015. Proceedings. Part
		II}, pages 171--183. Berlin: Springer, 2015.
	\newblock \href {https://doi.org/10.1007/978-3-662-47666-6_14}
	{\path{https://doi.org/10.1007/978-3-662-47666-6_14}}.
	
	\bibitem{reparam}
	U.~Fahrenberg and M.~Raussen.
	\newblock Reparametrizations of continuous paths.
	\newblock {\em J. Homotopy Relat. Struct.}, 2(2):93--117, 2007.
	
	\bibitem{zbMATH02231448}
	L.~Fajstrup.
	\newblock Dipaths and dihomotopies in a cubical complex.
	\newblock {\em Adv. Appl. Math.}, 35(2):188--206, 2005.
	\newblock \href {https://doi.org/10.1016/j.aam.2005.02.003}
	{\path{https://doi.org/10.1016/j.aam.2005.02.003}}.
	
	\bibitem{DAT_book}
	L.~Fajstrup, E.~Goubault, E.~Haucourt, S.~Mimram, and M.~Raussen.
	\newblock {\em Directed algebraic topology and concurrency. {With} a foreword
		by {Maurice} {Herlihy} and a preface by {Samuel} {Mimram}}.
	\newblock SpringerBriefs Appl. Sci. Technol. Springer, 2016.
	\newblock \href {https://doi.org/10.1007/978-3-319-15398-8}
	{\path{https://doi.org/10.1007/978-3-319-15398-8}}.
	
	\bibitem{FR}
	L.~Fajstrup and J.~Rosick{\'y}.
	\newblock A convenient category for directed homotopy.
	\newblock {\em Theory Appl. Categ.}, 21:7--20, 2008.
	
	\bibitem{model3}
	P.~Gaucher.
	\newblock A model category for the homotopy theory of concurrency.
	\newblock {\em Homology Homotopy Appl.}, 5(1):p.549--599, 2003.
	\newblock \href {https://doi.org/10.4310/hha.2003.v5.n1.a20}
	{\path{https://doi.org/10.4310/hha.2003.v5.n1.a20}}.
	
	\bibitem{mdtop}
	P.~Gaucher.
	\newblock Homotopical interpretation of globular complex by multipointed
	d-space.
	\newblock {\em Theory Appl. Categ.}, 22(22):588--621, 2009.
	
	\bibitem{266597}
	P.~Gaucher.
	\newblock About the category of all small diagrams.
	\newblock MathOverflow, 2017.
	\newblock URL:https://mathoverflow.net/q/266597 (version: 2017-04-07).
	\newblock \href {http://arxiv.org/abs/https://mathoverflow.net/q/266597}
	{\path{https://mathoverflow.net/q/266597}}.
	
	\bibitem{Moore1}
	P.~Gaucher.
	\newblock Homotopy theory of {M}oore flows ({I}).
	\newblock {\em {Compositionality}}, 3(3), 2021.
	\newblock \href {https://doi.org/10.32408/compositionality-3-3}
	{\path{https://doi.org/10.32408/compositionality-3-3}}.
	
	\bibitem{Moore2}
	P.~Gaucher.
	\newblock Homotopy theory of {M}oore flows ({II}).
	\newblock {\em {Extr. Math.}}, 36(2):157--239, 2021.
	\newblock \href {https://doi.org/10.17398/2605-5686.36.2.157}
	{\path{https://doi.org/10.17398/2605-5686.36.2.157}}.
	
	\bibitem{leftproperflow}
	P.~Gaucher.
	\newblock Left properness of flows.
	\newblock {\em Theory Appl. Categ.}, 37(19):562--612, 2021.
	
	\bibitem{QHMmodel}
	P.~Gaucher.
	\newblock Six model categories for directed homotopy.
	\newblock {\em Categ. Gen. Algebr. Struct. Appl.}, 15(1):145--181, 2021.
	\newblock \href {https://doi.org/10.52547/cgasa.15.1.145}
	{\path{https://doi.org/10.52547/cgasa.15.1.145}}.
	
	\bibitem{NaturalRealization}
	P.~Gaucher.
	\newblock Comparing cubical and globular directed paths.
	\newblock {\em {Fund. Math.}}, 262(3):259--286, 2023.
	\newblock \href {https://doi.org/10.4064/fm219-3-2023}
	{\path{https://doi.org/10.4064/fm219-3-2023}}.
	
	\bibitem{ThickCubes}
	P.~Gaucher.
	\newblock Towards a theory of natural directed paths, 2023.
	\newblock \href {https://doi.org/10.48550/arXiv.2306.02792}
	{\path{https://doi.org/10.48550/arXiv.2306.02792}}.
	
	\bibitem{DirectedDegeneracy}
	P.~Gaucher.
	\newblock Directed degeneracy maps for precubical sets.
	\newblock {\em Theory Appl. Categ.}, 41(7):194--237, 2024.
	
	\bibitem{Moore3}
	P.~Gaucher.
	\newblock Homotopy theory of {M}oore flows ({III}).
	\newblock {\em North-West. Eur. J. Math.}, (10):55--113, 2024.
	
	\bibitem{RegularMoore}
	P.~Gaucher.
	\newblock {Regular directed path and Moore flow}.
	\newblock {\em Rend. Mat. Appl., VII. Ser.}, 45(1-2):111--151, 2024.
	
	\bibitem{diCW}
	P.~Gaucher and E.~Goubault.
	\newblock Topological deformation of higher dimensional automata.
	\newblock {\em Homology Homotopy Appl.}, 5(2):39--82 (electronic), 2003.
	\newblock Algebraic topological methods in computer science (Stanford, CA,
	2001).
	\newblock \href {https://doi.org/10.4310/HHA.2003.v5.n2.a3}
	{\path{https://doi.org/10.4310/HHA.2003.v5.n2.a3}}.
	
	\bibitem{mg}
	M.~Grandis.
	\newblock Directed homotopy theory. {I}.
	\newblock {\em Cah. Topol. G\'eom. Diff\'er. Cat\'eg.}, 44(4):281--316, 2003.
	
	\bibitem{ref_model2}
	P.~S. Hirschhorn.
	\newblock {\em Model categories and their localizations}, volume~99 of {\em
		Mathematical Surveys and Monographs}.
	\newblock American Mathematical Society, Providence, RI, 2003.
	\newblock \href {https://doi.org/10.1090/surv/099}
	{\path{https://doi.org/10.1090/surv/099}}.
	
	\bibitem{MR99h:55031}
	M.~Hovey.
	\newblock {\em Model categories}.
	\newblock American Mathematical Society, Providence, RI, 1999.
	\newblock \href {https://doi.org/10.1090/surv/063}
	{\path{https://doi.org/10.1090/surv/063}}.
	
	\bibitem{MR2545830}
	S.~Krishnan.
	\newblock A convenient category of locally preordered spaces.
	\newblock {\em Appl. Categ. Structures}, 17(5):445--466, 2009.
	\newblock \href {https://doi.org/10.1007/s10485-008-9140-9}
	{\path{https://doi.org/10.1007/s10485-008-9140-9}}.
	
	\bibitem{MR1031717}
	M.~Makkai and R.~Par\'e.
	\newblock {\em Accessible categories: the foundations of categorical model
		theory}, volume 104 of {\em Contemporary Mathematics}.
	\newblock American Mathematical Society, Providence, RI, 1989.
	\newblock \href {https://doi.org/10.1090/conm/104}
	{\path{https://doi.org/10.1090/conm/104}}.
	
	\bibitem{376474}
	user168706.
	\newblock {Continuous bijection between two homotopy equivalent
		$\Delta$-generated spaces}.
	\newblock MathOverflow.
	\newblock URL:https://mathoverflow.net/q/376474 (version: 2020-11-14).
	\newblock \href {https://mathoverflow.net/q/376474}
	{\path{https://mathoverflow.net/q/376474}}.
	
	\bibitem{MR4070250}
	K.~Ziemia\'{n}ski.
	\newblock Spaces of directed paths on pre-cubical sets {II}.
	\newblock {\em J. Appl. Comput. Topol.}, 4(1):45--78, 2020.
	\newblock \href {https://doi.org/10.1007/s41468-019-00040-z}
	{\path{https://doi.org/10.1007/s41468-019-00040-z}}.
	
	\bibitem{Ziemiaski2012}
	K.~Ziemiański.
	\newblock Categories of directed spaces.
	\newblock {\em {Fund. Math.}}, 217(1):55–71, 2012.
	\newblock \href {https://doi.org/10.4064/fm217-1-5}
	{\path{https://doi.org/10.4064/fm217-1-5}}.
	
\end{thebibliography}

\end{document}